\documentstyle[11pt,amstex]{article}                                
\input epsf

\def\AFOUR{%
\setlength{\textheight}{9.0in}%
\setlength{\textwidth}{5.75in}%
\setlength{\topmargin}{-0.375in}%
\hoffset=-.5in%
\renewcommand{\baselinestretch}{1.17}%
\setlength{\parskip}{6pt plus 2pt}%
}

\AFOUR                                           

\makeatletter
\@addtoreset{equation}{section}

\makeatother

\newcommand{\nc}{\newcommand}
\newcommand{\rnc}{\renewcommand}

\nc{\bea}{\begin{eqnarray}}
\nc{\eea}{\end{eqnarray}}
\nc{\be}{\bea}
\nc{\ee}{\eea}
\rnc{\a}{\alpha}
\nc{\ab}{\bar{\a}}
\nc{\bs}{\backslash}
\nc{\ap}{\a^{+}}
\nc{\abm}{\ab^{-}}
\rnc{\b}{\beta}
\nc{\bb}{\bar{\b}}
\nc{\bbp}{\bb_{\zb}^{+}}
\nc{\bm}{\b_{z}^{-}}
\nc{\oa}{\overline{\a}}
\nc{\ob}{\overline{\b}}
\rnc{\gg}{\gamma}
\rnc{\d}{\delta}
\nc{\f}{\phi}
\nc{\fb}{\bar{\phi}}
\nc{\vf}{\varphi}
\nc{\p}{\psi}

\rnc{\c}{\chi}
\nc{\la}{\lambda}
\nc{\m}{\mu}
\nc{\n}{\nu}
\rnc{\o}{\omega}
\nc{\Om}{\Omega}
\rnc{\t}{\theta}
\nc{\eps}{\epsilon}
\rnc{\S}{\Sigma}
\nc{\F}{\Phi}

\nc{\trac}[2]{{\textstyle\frac{#1}{#2}}}

\nc{\ex}[1]{\mbox{e}^{\,\textstyle#1}}

\nc{\mat}[4]{\left(\begin{array}{cc}#1&#2\\#3&#4\end{array}\right)}

\nc{\som}[9]{\left(\begin{array}{ccc}#1&#2&#3\\#4&#5&#6\\#7&#8&#9%
\end{array}\right)}

\nc{\tr}{\mathop{\mbox{tr}}\nolimits}
\nc{\ad}{\mathop{\mbox{ad}}\nolimits}
\nc{\Tr}{\mathop{\mbox{Tr}}\nolimits}
\nc{\Det}{\mathop{\mbox{Det}}\nolimits}
\nc{\rk}{\mathop{\mbox{rk}}\nolimits}
\nc{\ra}{\rightarrow}
\nc{\Ra}{\Rightarrow}
\nc{\LRa}{\Leftrightarrow}
\nc{\ot}{\otimes}
\rnc{\ss}{\subset}
\nc{\nul}{\noindent\underline}
\nc{\non}{\nonumber\\}
\nc{\K}{K\"{a}hler}

\nc{\subs}[1]{{\vspace*{0.5cm}}%
{\noindent\underline{#1}}{\addcontentsline{toc}{subsection}{#1}}%
{\vspace*{0.3cm}}}

\nc{\zb}{\bar{z}}
\rnc{\lg}{\frak{g}}
\nc{\lt}{\frak{t}}
\nc{\lk}{\frak{k}}
\nc{\lh}{\frak{h}}
\nc{\pik}{\Pi_{\lk}}
\nc{\pip}{\Pi_{+}}
\nc{\pim}{\Pi_{-}}
\nc{\pih}{\Pi_{\lh}}
\nc{\jz}{J_{z}}
\nc{\jzh}{\jz^{\lh}}
\nc{\jzp}{\jz^{+}}
\nc{\jzm}{\jz^{-}}
\nc{\del}{\partial}
\nc{\dz}{\del_{z}}
\nc{\dzb}{\del_{\bar{z}}}
\nc{\az}{A_{z}}
\nc{\azb}{A_{\bar{z}}}
\nc{\g}{g^{-1}}
\nc{\dw}{\Delta_{W}}
\nc{\Ad}{{\mbox{Ad}}}
\nc{\ks}{Ka\-za\-ma-\-Su\-zu\-ki}
\nc{\KS}{\ks}
\nc{\ksm}{\ks\ model}

\nc{\BB}{{\Bbb B}}
\nc{\CC}{{\Bbb C}}
\nc{\PP}{{\Bbb P}}
\nc{\cpm}{\CC\PP(m)}
\nc{\cpn}{\CC\PP(n)}
\nc{\cp}[1]{\CC\PP(#1)}
\nc{\gmn}{G(m,m+n)}
\nc{\gmnk}{\gmn_{k}}
\nc{\cO}{{\cal O}}
\nc{\bcO}{\bar{\cO}}
\nc{\bO}{\bar{O}}
\nc{\oQ}{\overline{Q}}
\nc{\ZHS}{{\Bbb Z}{\mathrm HS}}
\nc{\QHS}{{\Bbb Q}{\mathrm HS}}
\nc{\Pf}{{\mathrm Pfaff}}

\catcode`\@=11
\long\def\@makecaption#1#2{%
    \vskip 10pt

\setbox\@tempboxa\hbox{
      \small\sf{\bfcaptionfont #1. }\ignorespaces #2}%
    \ifdim \wd\@tempboxa >\captionwidth {%
        \rightskip=\@captionmargin\leftskip=\@captionmargin
        \unhbox\@tempboxa\par}%
      \else
        \hbox to\hsize{\hfil\box\@tempboxa\hfil}%
    \fi
}
\font\bfcaptionfont=cmssbx10 scaled \magstephalf
\newdimen\@captionmargin\@captionmargin=2\parindent
\newdimen\captionwidth\captionwidth=\hsize
\catcode`\@=12

\newtheorem{theorem}{Theorem}[section]
\newtheorem{proposition}[theorem]{Proposition}
\newtheorem{lemma}[theorem]{Lemma}
\newtheorem{corollary}[theorem]{Corollary}

\newtheorem{definition}[theorem]{Definition}
\newtheorem{example}[theorem]{Example}
\newtheorem{remark}[theorem]{Remark}
\newtheorem{other}[theorem]{}

\nc{\bother}{\begin{other}}
\nc{\eother}{\end{other}}
\nc{\bex}{\begin{example}}
\nc{\eex}{\end{example}}
\nc{\bth}{\begin{theorem}}
\nc{\bdf}{\begin{definition}}
\nc{\edf}{\end{definition}}
\nc{\bpr}{\begin{proposition}}
\nc{\epr}{\end{proposition}}
\nc{\blm}{\begin{lemma}}
\nc{\elm}{\end{lemma}}
\nc{\br}{\begin{remark}}
\nc{\er}{\end{remark}}
\nc{\bcor}{\begin{corollary}}
\nc{\ecor}{\end{corollary}}

\def\lmn#1{\vadjust{\setbox1=\vtop{\hsize 12mm
\parindent=0pt\baselineskip=9pt
\rightskip=4mm plus 4mm#1}
\hbox{\kern-12mm\smash{\raise .5ex\box1}}}}

\def\K{\mathcal K}

\def\F{\mathcal F}

\def\bb{\beta}

\def\la{\langle}
\def\ra{\rangle}

\def\lbl#1{\label{#1}}

\def\noi{\noindent}

\def\f{\phi}

\def\t{\tau}

\def\d{\delta}
\def\p{\prime}

\def\Q{\mathbb Q}
\def\c{$\clubsuit$}

\def\I{{\mathrm I}}
\def\noi{\noindent}

\def\s{\sigma}

\def\build#1_#2^#3{\mathrel{\mathop{\kern 0pt#1}\limits_{#2}^{#3}}}

\begin{document}
\makeatletter
\begin{titlepage}
\begin{flushright}
IC/99/66
\end{flushright}
\begin{center}
{\LARGE\bf The
Universal Perturbative Quantum 3-manifold Invariant, Rozansky-Witten
Invariants, and the Generalized Casson Invariant.}

\vskip 0.4in
{\bf Nathan Habegger and George Thompson}
\vskip .1in

\end{center}

\begin{abstract}Let $Z^{LMO}$ be the 3-manifold invariant of \cite{LMO}. It is
shown that $Z^{LMO}(M)=1$, if the first Betti number of $M$, $b_1(M)$, is
greater than 3.  If $b_1(M)=3$, then $Z^{LMO}(M)$ is completely determined
by the cohomology ring of $M$.
A relation of $Z^{LMO}$ with the Rozansky-Witten invariants
$Z_{X}^{RW}[M]$ is established at a physical level of
rigour. We show that $Z_{X}^{RW}[M]$ satisfies appropriate
connected sum properties suggesting that the generalized Casson
invariant ought to be computable from the LMO invariant.
\end{abstract}

\end{titlepage}

\tableofcontents

\makeatother

\section{Introduction}

In [W], E. Witten explained the Jones polynomial using physics.
In doing so, he introduced mathematicians to the partition function
$Z_{G,k}^{CS}(M,L)$ of the topological quantum field theory
associated to the Chern Simons action, for a Lie group $G$ and
coloured\footnote{By coloured link, one means that to each link
component, one associates a representation of $G$} link $L\subset M$.
Its physical definition is given by a Feynman path integral
over the infinite dimensional space of connections.\footnote{The
connections are on
an underlying principle $G$-bundle lying over $M$.}

In general, one expects that topological field theories defined using
the path integral, or perturbative versions of these, can be given a
rigorous
definition through surgery formulae, just as is the case for other
quantum invariants, such as  the Reshetikhin-Turaev \cite{RT} invariants,
$Z^{RT}_{G,k}$, or the more recent universal invariant
$Z^{LMO}$, of T. Le, J. Murakami, and T. Ohtsuki \cite{LMO} and the
\AA{}rhus invariant, $Z^{\AA}$, \cite{BGRT}. Invariants
have also been given through integral formulae, as is the case for the
Kontsevich integral $Z^K$ (see \cite{K} \cite{B}), the Bott-Taubes
invariant $Z^{BT}$, \cite{BoTa}, \cite{AF}, and the
invariant $Z^{BC}$ of Bott and Cattaneo \cite{BC}, \cite{AC}. (It is
conjectured that $Z^{K}=Z^{BT}$.)

Our intent in this paper is to study the invariant
$Z^{LMO}(M)$. $Z^{LMO}(M)$ lies in $A(\emptyset)$, the vector
space of Feynman diagrams modulo anti-symmetry and $IHX$
relations. This vector space is not well-understood, except in low degrees.
(See Vogel, \cite{V}, for an attempt to understand the structure of
$A(\emptyset)$.)

Quantum invariants (or perturbative versions of these) are a rich
source of data for the study of knots, links, and
3-manifolds. Nevertheless, their relationship to classical topology
remains obscure, hampering their use in problem-solving. A notable
exception is the Alexander polynomial of a knot, which, through its
interpretation as the Conway polynomial (together with the solution of the
Conway weight system on uni-trivalent graphs \cite{KK}), gives a
computation of
(the one loop) part of the Kontsevich integral. Another recent advance
has been the computation of the Milnor invariants from the Kontsevich
integral \cite{HM}.

For 3-manifolds, one has the result that the degree one term of the LMO
invariant, $Z^{LMO}_{1}$, is the Casson-Walker-Lescop
invariant \cite{LMO}, \cite{BeHa}. However, beyond this, the
topological significance of $Z^{LMO}$ remains a mystery. For example, it is not
even known whether or not the degree two term of $Z^{LMO}$ vanishes in the
simply connected case (which of course would be implied by a positive
solution of the Poincar\'{e} conjecture, since $Z^{LMO}(S^{3})= 1$). 

One possible programme for attempting to tie the quantum invariants to
homotopy data is through generalization of the Casson invariant to
groups other than $SU(2)$, e.g., $SU(N)$. Recent advances on the
mathematical side \cite{BH}, for $SU(3)$, as well as on the physical side,
by Rozansky and Witten \cite{RW}, may make this programme
tractable. The purpose of this paper is to give a conjectural
relationship between the generalized
Casson invariants and $Z^{LMO}$ and some partial evidence for its
veracity. We consider this conjecture to be an explicit form of the
basic philosophical viewpoint of \cite{RW}, who believe their invariants
are of finite type.

Indeed, we may summarize the underlying ideas of \cite{RW} as
follows. On the one hand a comparison of the Rozansky-Witten invariants
to the perturbative Chern-Simons theory indicates that, for
$b_{1}(M)=0$, they both arise from one universal theory. The
difference between the two rests in the choice of weight system (in
\cite{RW} a rigorous mathematical weight system $W_{X}^{RW}$ is
given). On the other hand, (and perhaps the deepest part of the
theory) equivalence between certain physical theories allows one to identify
the Rozansky-Witten invariants for a particular choice of
hyper-K\"{a}hler manifold $X_{G}$ with a regularized Euler
characteristic \cite{BT1} of the moduli
space of flat G-connections on the 3-manifold M,
$\chi_{G}(M)$. The equivalence comes from the work of Seiberg and
Witten on 3-dimensional theories (and is the analogue of their
work in four dimensions). The twist of the first theory yields the
gauge theoretic model of the Euler characteristic while the twist of the
second is the Rozansky-Witten model. The equivalence of the physical
theories suggests the equivalence of their twisted topological
versions. Thus, from the physics side, one expects that
\be
Z^{RW}_{X_{G}}(M) = \chi_{G}(M)  . \label{eul}
\ee
One important consequence of (\ref{eul}) is its potential use in
computing $\chi_{G}(M)$.

On the mathematical side there are, currently, a number of candidates
for a universal perturbative quantum invariant. These include the LMO
invariant, $Z^{LMO}$, the \AA{}rhus invariant, $Z^{\AA}$, and the
invariant of Bott and Cattaneo, $Z^{BC}$. It has been suggested
\cite{BGRT} that although the LMO
invariant agrees with the \AA{}rhus invariant, for rational homology
spheres, that nevertheless $Z^{LMO}$ is more directly related to
the Rozansky-Witten sigma model theory than to the Chern-Simons gauge
theory.

For $b_{1}(M) >0$, (except for $b_{1}(M)=1$ and ${\mathrm
Tor\, 
H}_{1}(M, {\Bbb Z}) \neq 0$) the first author and collaborators,
\cite{BeHa} \cite{GH}, have
calculated $Z^{LMO}$ from classical data. Here, we perform analogous
computations
in the Rozansky-Witten theory and observe, for $b_{1}(M)>0$, that
these results agree. Specifically, we show, for $b_{1}(M) >0$, and under
the conditions mentioned above,\footnote{The computations that
we make for the Rozansky-Witten theory suggest that it is to be
expected that the results of
\cite{GH} hold even when the manifold has torsion in ${\mathrm
H}_{1}(M, {\Bbb Z})$.} that at the physical level of rigour,
\be
W^{RW}_{X}\left( Z_{n}^{LMO} (M) \right) = Z^{RW}_{X}(M) , \label{wz}
\ee
where $X$ denotes a hyper-K\"{a}hler manifold of dimension $4n$.

One might naively conjecture that (\ref{wz}) holds for all
3-manifolds. However, for $b_{1}(M)=0$, which is the case of most interest,
numerous considerations, including connected sum formulae and
normalization conditions, indicate that the equality (\ref{wz}) should be
modified. We introduce invariants $\lambda^{k}_{X}(M)$ for all $k$ and
$X$ which are computed from the Rozansky-Witten theory, and the
formulae now suggest\footnote{This corresponds to the fact that
$Z_{G,k}^{RT}/Z_{G,k}^{RT}(S^{3})$ and $\sum_{n} |H_{1}|^{-n} Z_{n}^{LMO}$
are both multiplicative.} that
\be
\left| {\mathrm H}_{1}(M, {\Bbb Z})\right|^{n-k}W^{RW}_{X}\left
( Z_{k}^{LMO} (M) \right) = \lambda_{X}^{k}(M) . \label{wl}
\ee

We now propose, that on correcting for the trivial connection, one
should replace (\ref{eul}) with the equality
\be
\lambda^{n}_{X_{G}}(M) = \lambda_{G}(M) ,
\ee
where $\lambda_{G}(M)$ is the, still to be mathematically defined, G-Casson
invariant, and where ${\mathrm rank}(G)=n$. In this way, we obtain the
purely mathematical 

{\bf Conjecture:}
\be
W^{RW}_{X_{G}}\left( Z_{n}^{LMO} (M) \right) = \lambda_{G}(M). \label{conjgc}
\ee
(N.b., for $SU(2)$ this equality holds by the computation of $Z^{LMO}$
combined with those on the physics side \cite{RW}.)

The equalities above are certainly suggestive.
On the one hand, $Z^{LMO}(M)$ satisfies axioms\footnote{Actually, the TQFT
axioms hold for certain truncations of $Z^{LMO}$.} of topological quantum
field theory (TQFT), see \cite{MO}, as is the case for the Chern-Simons
theory $Z^{RT}_{G,k}$. On the other hand,
$Z^{RW}_X$ is  given as a topological
sigma-model. As explained in \cite{RW}, the actions of the
Chern-Simons theory, and the Rozansky-Witten theory are formally
analogous (see section 4 below).

We begin this paper with a computation, which
originally appeared in
\cite{H1},  of $Z^{LMO}(M)$ for manifolds whose first Betti
number, $b_1(M)$, is greater than or equal to 3. Subsequently, computations
for $b_1(M)=2$, \cite{BeHa}, and $b_1(M)=1 $, \cite{GH},
followed. These computations were inspired by the
work of T. Le, \cite{Le}, who showed that the invariant $Z^{LMO}(M)$,
restricted to homology spheres, is the universal finite type
invariant\footnote{See \cite{H2} for an expository account of the theory of
finite type invariants.} in the sense of Ohtsuki \cite{O}. 

Specifically, in section 3, we will give a proof of parts $(i)$ and
$(ii)$ of the
following result (part $(iii)$ was proven in \cite{BeHa} and part $(iv)$ in
\cite{GH}). Let $\lambda_M$ denote the Lescop invariant of $M$, see \cite{L}.

\vskip .8cm
\noindent \bf Theorem 1. \it

(i) Suppose $b_1(M)>3$. Then $Z^{LMO}(M)=1$.

(ii) There are non-zero $\gamma_n\in A_n(\emptyset)$, such that if
$b_1(M)=3$, 
then $Z^{LMO}(M)=\Sigma_n \lambda_M^n \gamma_n$.

(iii) \cite{BeHa} There are non-zero ${\mathcal H}_n\in
A_{n}(\emptyset)$, such
that if $b_1(M)=2$, then \linebreak
$Z^{LMO}(M)=\Sigma_n \lambda_M^n {\mathcal H}_n$.

(iv) \cite{GH} For $H_1(M)={\Bbb Z}$, $Z^{LMO}(M)$ determines and is
determined by $A(M)$, the Alexander Polynomial of $M$.

\rm
\noindent{\bf Remark.} In fact, though not observed in \cite{BeHa}, but
as suggested from the combinatorics of the  physical approach, one can
show using equality (2) in
\cite{BeHa} that $\gamma_{n} =\pm {\mathcal H}_{n}$.

Sections \ref{sec.theory}-\ref{manb1} of this paper concern heuristic
results, reminiscent of
theorem 1 and the hypothetical equality
$W_X^{RW}(Z^{LMO}_n(M))=Z_X^{RW}(M)$.  Specifically, we give a 
heuristic proof, (i.e. at the physical level of rigour) of the following:

\vskip .8cm
\noindent \bf Heuristic Theorem 2. \it

(i)  Suppose $b_1(M)>3$.  Then $Z_X^{RW}(M)=0$.

(ii)  There are constants $c_X$, such that if  $b_1(M)=3$, then
$Z_X^{RW}(M)=c_X \lambda_M^n$.

(iii)  There are constants $c'_X$, such that if  $b_1(M)=2$, then
$Z_X^{RW}(M)=c'_X \lambda_M^n$.

(iv)  For $b_1(M)=1$, $Z_X^{RW}(M)$ is determined by Reidemeister Torsion.
\rm
\vskip .8cm

Actually, the equality $W_X^{RW}(Z^{LMO}_n(M))=Z_X^{RW}(M)$
suggests that $c_X=W_X^{RW}(\gamma_n)$ and $c'_X=W_X^{RW}({\mathcal
H}_n)$. Our calculations indicate that this is so and furthermore,
show that $c_{X}=c_{X}'$.

The final sections, \ref{manb10}-\ref{gcasson}, are devoted to
deriving the properties of the $\lambda^{k}_{X}$ invariants that are
required to motivate (\ref{wl}) and our conjecture (\ref{conjgc}).

Let us conclude this introduction with a few remarks on our proof of
heuristic theorem 2. While, for some purposes, the perturbative Feynman
diagram expansion may be useful, e.g., for obtaining weight systems,
our approach is essentially non-perturbative. In general the path
integral formalism may have uses beyond giving us a definition of
invariants. One can define theories via path integrals and after
passing to the perturbation theory completely forgo the path integral
formulation. This leads to an
interesting set of combinatorial problems, having to do with the type
of diagrams to be considered, as well as their frequency. On
the other hand, it may happen that the path integral can be performed
in a, more or less, elementary manner. In this case the
combinatorial issues are by-passed, and in addition one obtains nicely
re-summed formulae. An example of such a situation is the derivation of
the Verlinde formula \cite{BT2} for the dimension of the space of
holomorphic sections of the
$k$'th tensor power of the determinant line bundle over the space of
flat connections on a Riemann surface. Similarly, for the
Rozansky-Witten invariants, we will see that it is better to `perform' the path
integral directly rather than to expand out first.

The main thrust of our physical computations is then to avoid
working directly with diagrams. However, in order to make the
relationship with
\cite{LMO} somewhat more transparent we will, on occasion, explain
certain phenomena at the diagrammatic level.

\noindent
{\small
\bf Acknowledgments:\rm\ \  This paper was begun at the
Mittag-Leffler Institute, where the authors were participants in the
special year on topology and physics.  We thank the institute for its
support, and its staff, for their friendly and efficient professional
assistance. Thanks are also due to the ICTP for support. 
We also extend thanks to J.~Andersen, M.~Blau, H.~Murakami,
D.~Pickrell, and S.~Rajeev,
for the stimulating conversations we had with them during the elaboration
of this paper. This work was
supported in part by the EC under the TMR contract ERBF MRX-CT 96-0090.}

\section{The Invariant $Z^{LMO}(M)$.} \lbl{sec.invariant}

The invariant $Z^{LMO}(M)$ is computed in general from the Kontsevich integral
(denoted here by $Z^{K}(L)$) of any framed link $L\colon
\coprod_{j=1}^{j=l}S^1\rightarrow {\bf R}^3$, such that surgery on $L$,
denoted by $S^3(L)$, produces $M$.  $Z^{K}(L)$ lies in
$A(\coprod_{j=1}^{j=l}S^1)$.

Before stating the result, we recall (see, e.g., \cite{B}, \cite{LM1},
\cite{LM2}, \cite{V}) that
$A(X)$ denotes the graded-completed $\Q$-vector space of Feynman
diagrams
$X\cup \Gamma$ on the compact 1-manifold $X$.  The space $A(X)$ is
graded by the degree, where the degree of a diagram is half the number
of vertices of $\Gamma$.   Using the notation of [HM], we let
$\coprod_{j=1}^{j=l} \I_j$ denote the disjoint union of $l$ copies of the
interval, and we set $A(l)= A(\coprod_{j=1}^{j=l} \I_j)$.   $A(l)$ is a Hopf
algebra, and one has that $A(1)=A(S^1)$.  Moreover,  any embedding
$I\rightarrow X$ gives rise to a well defined action of $A(1)$ on $A(X)$.  In
particular, $A(1)^{\otimes l}$ acts on $A(l)$ and on
$A(\coprod_{j=1}^{j=l}S^1)$.

For $1\le i,j \le l$, we let $\xi_{ij}$ be the degree 1 diagram
$X\cup \Gamma$, with $X=\coprod_{j=1}^{j=l} \I_j$, where $\Gamma$ is a
chord with vertices on the $i$-th and $j$-th components ($i$ may be
equal to $j$).  We set $\xi_{123}=[\xi_{12},\xi_{23}]$, where $[a,b]=ab-ba$
denotes the Lie bracket of $a$ and $b$.  ($\xi_{123}$ is represented by the
diagram $X\cup \Gamma$, with $X=\coprod_{j=1}^{j=3} \I_j$,
where $\Gamma$ is the $Y$-graph of degree 2 with one vertex on each
component of $X$.)

In \cite{LMO}, maps $i_n\colon A_{nl+i}(\coprod_{j=1}^{j=l}S^1)\rightarrow
A_i(\emptyset)$ were defined for $i\le n$.    (We set $i_n=0$ otherwise.)
We denote by
$p_l\colon A(l)\rightarrow A(\coprod_{j=1}^{j=l}S^1)$ the quotient
mapping. We set $\gamma_0=1$, and $\gamma_n=i_n
( p_3({ {\xi_{123}^{2n} } \over {(2n)!} }))\in A_n(\emptyset)$.
Note that $A_1(\emptyset)$ is 1-dimensional and
$A_1(\emptyset)^{\otimes n}$ is a direct summand of $A_n(\emptyset)$.
Moreover, it is easily seen from the definition of $i_n$ that the image of
$\gamma_n$ in $A_1(\emptyset)^{\otimes n}$ is nonzero.  Hence
$\gamma_n$ is nonzero.

For a set $A$, we set $|A|$ to be the cardinality of
$A$, if this is finite, and $0$ otherwise.  For a 3-manifold $M$ with
$b_1(M)=3$, we define
$\lambda_M=|Tor(H_1(M))|\ |{{H^3(M)}\over{i(H^1(M)^{\otimes 3})}}|^2$,
where $i\colon H^1(M)^{\otimes 3}\rightarrow H^3(M)$ is given by the
cup product $a\otimes b \otimes c\mapsto a\cup b \cup c$.  This is
Lescop's
invariant, for $b_1(M)=3$, see \cite{L} section 5.3.

\vskip .8cm 
\noindent \bf{Theorem 1}. \it

(i)  Suppose $b_1(M)>3$.  Then $Z^{LMO}(M)=1$.

(ii)  Suppose $b_1(M)=3$.  Then $Z^{LMO}(M)=\Sigma_n \lambda_M^n \gamma_n$.
\rm
\vskip .5cm

The theorem will be proven in the next section. We first recall here how
$Z^{LMO}(M)$ is defined.  One puts $\nu=Z^{K}(U)$, where $U$ is the trivial knot with
faming zero. Then $Z^{LMO}_n(M)$ is the degree $n$ part of the expression

\be
{i_n(Z^{K}(L)\nu^{\otimes l})}\over {(i_n(\nu^2 exp({{\xi_{11}}\over
2}))^{b_+} (i_n(\nu^2 exp({-{\xi_{11}}\over 2}))^{b_-}} \label{om}
\ee
considered to lie
in  $A_{\le n}(\emptyset)$, where $b_+, b_-$ denote the number of positive
and negative eigenvalues of the linking matrix of $L$.  (It was shown in
\cite{LMO} that the expressions in the denominators are invertible.) 

\noindent{\bf Remark.} For later use we note that for $k
\leq n$, in \cite{LMO} the degree k part of the above expression
(\ref{om}) was denoted by
$\Omega_{n}(M)^{(k)}$. Thus in particular, $Z_{n}^{LMO}(M) =
\Omega_{n}(M)^{(n)}$.

In the proof of the theorem, we will need to make use of certain facts.

1)  $i_n$ satisfies the property that it vanishes on diagrams which have
fewer than $2n$ vertices on some component.

2)  Let $\s$ be a string link whose closure is $L$.  Then
$Z^{K}(L)=p_l(Z^{K}(\s)\nu_l)$ (see \cite{LM2}). Here $Z^{K}(\s)$ lies
in $A(l)$ and $\nu_l\in
A(l)$ is obtained from $\nu=\nu_1$ by the operator which takes a diagram
on the interval  to the sum of all lifts of vertices to each of the $l$
intervals.  It is known that $\nu$ and hence $\nu_l$ is a sum of diagrams
each of which has each component of $\Gamma$ non-simply connected
(see \cite{HM}).

3)  Let $\s$ be an $l$-component string link.  Then
$Z^{K}(\s)=exp (\xi^t +\xi^h)$, where $\xi^t$ is a linear combination of
diagrams for which $\Gamma$ is a tree, and $\xi^h$ is a linear
combination of diagrams for which $\Gamma$ is connected, but not
simply connected.  If we denote by $A^t(X)$ the quotient of $A(X)$
obtained by setting to zero all diagrams for which some component of
$\Gamma$ is not simply connected, and denote by $Z^t(\s)$, the image in
$A^t(l)$ of $Z^{K}(\s)$, then it was shown in \cite{HM} that the
Milnor invariants of
$\s$ determine, and are determined by $Z^t(\s)=exp (\xi^t )$. We will
need the fact that if the linking numbers and framings are zero, then
$\xi^t$ has degree $\ge 2$, and moreover, the coefficient of
$\xi_{123}$ is the Milnor invariant $\mu_{123}$.  (See \cite{HM}).

\section{Proof of Theorem 1.} \lbl{sec.proof1}

The theorem will be proven progressively, starting from the case where
$M$ is obtained via surgery on  an algebraically split link $L$ (i.e., one with
vanishing linking numbers) having 3 components all of
which are zero-framed.  In this case, $H_1(M)={\Bbb Z}^3$, so that
$Tor(H_1(M))=0$. Moreover, using the Poincar\'e dual interpretation of
cup product, one easily checks from the definition of the Milnor invariant
$\mu=\overline \mu_{123}(L)$ in terms of intersections of Seifert
surfaces (which can be completed to surfaces in $M$), that $|\mu|=
|{{H^3(M)}\over{i(H^1(M)^{\otimes 3})}}|$.  It follow that
$\lambda_M=\mu^2$ in this case.

The theorem in this case is an immediate
consequence of the observation that the only term contributing to
$Z^{LMO}_n(M)$ is $p_3({{\mu^{2n}\xi_{123}^{2n}}\over{(2n)!}})$. To
see this let
$\s$ be a zero framed string link whose closure is $L$.  Then by 3) above
and \cite{HM} (since the linking and framings are zero and $\mu=\overline
\mu_{123}(L)$),  $Z^{K}(\s)=exp(\mu\xi_{123}+\xi')$, where $\xi'$ is a
linear combination of diagrams, all of which consist of diagrams for which
$\Gamma$ is either not simply connected, or is a tree of degree $\ge 3$.
Note that in each case, such a diagram has a ratio of external vertices
(the univalent vertices of $\Gamma$) to internal vertices which is $<3$,
whereas this ratio for $\xi_{123}$ is $3$.  It follows that every term of
$Z^{K}(L)\nu^{\otimes l}=p_3(Z^{K}(\s)\nu_3)\nu^{\otimes l}$, which
has at least
$2n$ vertices on every component, must have at least ${{6n}\over 3}=2n$
internal vertices.  Hence such a term has degree at least ${{6n+2n}\over
2}=4n$, and has degree precisely
$4n$ if and only if that term is
$p_3({{\mu^{2n}\xi_{123}^{2n}}\over{(2n)!}})$.

Now suppose that $M=S^3(L)$, where $L$ is an algebraically split link
having $l$ components, and such that
$L$ contains a 3-component sublink $L_0$ which is zero-framed.  We set
$L_1=L\setminus L_0$.

We first suppose that $L_0$ and $L_1$ are separated by a 2-sphere, so
that $M$ is a connected sum. Recall from (\cite{LMO}, 5.1), that if $M$ is a
connected sum of $M'$ and $M''$ such that $b_1(M')>0$, then one has the
formula $Z^{LMO}_n(M)=Z^{LMO}_n(M')|H_1(M'')|^n$.  Setting $M'=S^3(L_0)$ and
$M''=S^3(L_1)$, then this formula shows that the result for $M$
is implied by the result for $M'$
shown earlier. (This includes the vanishing result if $b_1(M)>3$, since in
this case $b_1(M')>0$, and hence $|H_1(M'')|=0$.)

If  $L_0$ and $L_1$ are not separated by a 2-sphere, i.e., $L\ne L_0\coprod
L_1$, the result still follows, since one has that $i_n(Z^{K}(L)\nu^{\otimes l})=
i_n(Z^{K}( L_0\coprod L_1)\nu^{\otimes l})$.  To see this, let $\s$ be a string
link, whose closure is $L$, such that the first $3$ components,
$\s_0$, close up to give $L_0$.  Set $\s_1=\s\setminus \s_0$.   Let
$\s_0\times\s_1$ denote the juxtaposition of $\s_0$ and $\s_1$.
We wish to compare $Z^{K}(\s)$ and $Z^{K}(\s_0\times\s_1)$.  One has that
$Z^{K}(\s_0)=exp(\xi_0)$, $Z^{K}(\s_1)=exp(\xi_1)$, and hence that
$Z^{K}(\s_0\times\s_1)=exp(\xi_0+\xi_1)$.  Moreover, $Z^{K}(\s)=
exp(\xi_0+\xi_1+\xi') $, where $\xi'$ is a sum of diagrams for which
$\Gamma$ is connected, has degree $\ge 2$ (since $L$ is algebraically split),
and has a vertex on $\s_0$ and on $\s_1$.  Since $\xi_0$ is also of degree
$\ge 2$, it follows that every term of $(Z^{K}(L)-Z^{K}( L_0\coprod L_1))
\nu^{\otimes l}=p_l((Z^{K}(\s)-Z^{K}(\s_0\times\s_1))\nu_l)\nu^{\otimes l}$
is a sum of terms which satisfy that each component of $\Gamma$, with a
vertex on one of the 3 components of $L_0$, has degree at least 2 and that
some such component must also have a vertex lying on $L_1$.  It follows
that any such term, having at least $2n$ external vertices on each
component of $L$, must have more than $2n$ internal vertices, and
hence that such a term is in the kernel of $i_n$ (since it is of degree
$>nl+n$).

Now suppose that $M=S^3(L)$, where $L$ is arbitrary.  Let $B$ be the
linking matrix of $L$. It is well known that $B$ becomes diagonalizable after
taking the direct sum with a certain diagonal matrix $D$ having non-zero
determinant. Let $L'$ denote a link whose linking matrix is $D$.   Then if
$M''$ denotes $S^3(L\coprod L')$, the theorem holds for $M''$, since
$L\coprod L'$ is equivalent to an algebraically split link $L''$, via handle
sliding (so that $M''=S^3(L'')$). Then the theorem holds also for $M$, using
the formula $Z^{LMO}_n(M'')=Z^{LMO}_n(M)|H_1(S^3(L')|^n$ (since
$|H_1(S^3(L')|=|det(B)|\ne0$).

\section{Review of Rozansky-Witten Theory.} \label{sec.theory}

The theory whose partition function is believed to yield the
$G$-Casson invariant is a twisted version of $N=4$ super-Yang-Mills
theory in 3-dimensions \cite{W1}, \cite{BT1}. Seiberg and Witten
\cite{SW} have given a solution of the physical theory with $G=SU(2)$
in the coulomb branch. The coulomb branch of a theory corresponds to
an analysis at a particular (low) energy scale. This solution has, as its
moduli space, the reduced $SU(2)$ 2-Monopole moduli space, that is the
Atiyah-Hitchin space $X_{AH}$. Since the topological theory should not
depend on which scale we are looking at, we can twist the low energy
theory of Seiberg and Witten and in this way we are led to equating the $SU(2)$
Casson invariant with a particular path integral over the space of
maps from a 3-manifold to $X_{AH}$. 

Rather more generally it is
believed that the moduli space for for the physical theory with group
$G$ is some monopole moduli space. For example for $SU(n)$ it is
believed to be the reduced $SU(2)$ n-monopole moduli space. These
moduli spaces are all hyper-K\"{a}hler. We denote those
hyper-K\"{a}hler manifolds that arise as the moduli space of the coulomb
branch  of the $G$ physical theory by $X_{G}$. From this point of view the
$G$ Casson invariant is then equated with a particular path integral
over the space of maps from a 3-manifold to some hyper-K\"{a}hler
$X_{G}$. The path integral in question, $Z^{RW}_{X_{G}}[M]$, was described and
analysed in \cite{RW}. Given some subtleties that we will address
later, one expects that $Z^{RW}_{X_{G}}[M]$ and $\lambda_{G}(M)$ if not
equal are very closely related. (The exact statement was given in the
introduction (\ref{wl}).)

\subsection{The Rozansky-Witten Model}

Rozansky and Witten \cite{RW} defined a path integral, and so
invariants for a 3-manifold, for
any hyper-K\"{a}hler $X$. This section is devoted to describing the
objects that go into defining that path integral.

Let $\f$ be a map from the 3-manifold $M$ to a hyper-K\"{a}hler
manifold $X$. In local coordinates on $X$, the map is denoted
$\f^{i}$, $i= 1, \dots , 4n$.\footnote{$\f^{i}$ is the composite of
$\f$, restricted to
the inverse image of the coordinate neighborhood, with
the $i$-th coordinate function. Thus $\f^{i}$ is not a function defined on
all of $M$, but only on some open set.} We write $\f \in {\mathrm
Map}(M,X)$. Denote by ${\mathrm T_{\f}\, Map}(M,X)$ the tangent space of
${\mathrm Map}(M, X)$ at $\f$. One may identify $T_{\f}\,
{\mathrm Map}(M,X)$ with $\Omega^{0}\left(M, \f^{*}\left( TX
\right)\right)$.
Since one has, $TX \equiv_{{\Bbb R}}
T^{(1,0)}X \equiv V $, for  
the tangent bundle of $X$ (see the review in Appendix \ref{hkprops}). We
define $\eta$ to be a Grassman variable\footnote{The definition of
what it means to be a Grassman variable on a vector space
is explained in Appendix~\ref{sec.manipulations}} on $\Omega^{0}\left( M,
\f^{*}\left( V \right)\right)$, that is $\eta \in \Lambda^{1}
\left(\Omega^{0}\left(M, \f^{*}\left( V\right)\right)^{*}  \right)$
which is, in local coordinates, denoted
by $\eta^{I}(x)$, $I= 1, \dots, 2n$.  Let $\chi$ be a Grassman variable on
$\Omega^{1}\left(M ,\f^{*}\left(V\right)\right)$, that is it is an
element of $\Lambda^{1} \left(\Omega^{1}\left(M, \f^{*}\left( V
\right) \right)^{*} \right)$ which, in local coordinates, we denote by
$\chi^{I}$.

We define a Lagrangian (density on $M$) $L={\mathbf L}_1 +{\mathbf L}_2$,
\bea
{\mathbf L}_{1} & = & \frac{1}{2}g_{ij}(\f) \,d\f^{i}*d\f^{j}
+ \eps_{IJ} (\f) \, \chi^{I}* \nabla \eta^{J} \label{lag1}\\
{\mathbf L}_{2} & = & \frac{1}{2}\left( \eps_{IJ}(\f) \, \chi^{I}
\nabla \chi^{J} +
\frac{1}{3} \Omega_{IJKL}(\f) \, \chi^{I}\chi^{J}\chi^{K} \eta^{L}
\right). \label{lag2} 
\eea
The covariant derivative $\nabla$ is defined with the pullback of the
Levi-Civita connection on $V$,
\be
\nabla_{\, J}^{I} = d\d^{I}_{J} +
(d\f^{i})\Gamma^{I}_{i J} ,
\ee
and $*$ is the Hodge star operator on $M$ thought of as a Riemannian
manifold. 
The two Lagrangians are separately invariant under a pair of BRST
transformations. One does not need to pick a complex structure to
exhibit these, however that level of generality is not required and we
pick now a complex structure on $X$ so that the $\f^{I}$ are local
holomorphic coordinates with respect to this complex structure. In
this complex structure we can pick a basis, $Q$, $\overline{Q}$ for
the BRST charges which act by
\be
\begin{array}{ll}
\overline{Q}\f^{I} = 0, & \overline{Q}\f^{\overline{I}} =
T^{\overline{I}}_{\; J}\, \eta^{J}, \\
\overline{Q}\eta^{I}= 0, & \overline{Q}\chi^{I} = -d\f^{I} ,
\end{array} \label{brst1}
\ee
and
\be
\begin{array}{ll}
Q\f^{I} = \eta^{I}, & Q\f^{\overline{I}} = 0, \\
Q\eta^{I}= 0, & Q\chi^{I} =
T^{I}_{\overline{J}}\, d\f^{\overline{J} } - \Gamma_{JK}^{I}\,\eta^{J}\chi^{K} .
\end{array} \label{brst2}
\ee
These BRST operators satisfy the algebra
\be
Q^{2}= 0, \;\;\; \{ Q\, ,\, \overline{Q} \} = 0, \;\;\; \overline{Q}^{2} =
0.
\ee
The BRST invariant sigma model action is
\be
S = \int_{M} \left( {\mathbf L}_{1} + {\mathbf L}_{2} \right)
. \label{sigmod} 
\ee

We note that ${\mathbf L}_{1}$ is both $Q$ and $\overline{Q}$
exact. Indeed one has
\bea
{\mathbf L}_{1} &=& \langle d\f , d\overline{\f} \rangle +\langle \chi
, T \nabla \eta \rangle \nonumber \\
 & = & -\overline{Q} \langle \chi , d\overline{\f} \rangle \nonumber \\
& =&  \overline{Q} \, Q\; \left(\frac{1}{2}\eps_{IJ} \,
\chi^{I}  * \chi^{J}  \right) , \label{dbrst}
\eea
where the inner product for $X \in \Omega^{1}(M, T^{(1,0)}X )$ and $Y \in
\Omega^{1}(M,  T^{(0,1)}X  )$, is defined to be 
\be
\langle X , Y \rangle = g_{I \overline{J}} \, X^{I} * Y^{\overline{J}} .
\ee
In order to write ${\mathbf L}_{1}$ we needed to pick a metric on
$M$. However, as ${\mathbf L}_{1}$ is BRST exact, nothing depends on
the choice made (see \cite{BBRT} section 2 for how this is
established) and, ultimately, this explains why this theory produces a
3-manifold invariant.

Now we have a gauge theory interpretation of the sigma model action
(\ref{sigmod}) as a gauge fixed action. Firstly, ${\mathbf L}_{2}$, is
BRST invariant (and metric independent). However, it is not BRST
exact. So we may consider it to be the initial gauge invariant
Lagrangian that needs to be augmented with a gauge fixing term, in
order to arrive at a well defined theory. The
gauge fixing term should be BRST exact and we see that ${\mathbf
L}_{1}$ fits the bill. In section \ref{pirs} we will, for
$b_{1}(M)=1$, take this point of view and gauge fix an invariant
Lagrangian. In fact what one finds is a theory that looks a great deal
like the Chern-Simons theory of Witten. This suggestive analogy will
be taken up again when we make a more comprehensive comparison with
Chern-Simons theory below.

The action (\ref{sigmod}) at first sight defines quite a complicated
theory. However, as it is a topological theory, one may expect rather
drastic simplifications. This is, indeed, the case.

There are various
arguments that are available (see \cite{RW} and
\cite{T}) that establish that one may as well instead consider the
Lagrangians 
\bea
{\mathbf L}_{1} \rightarrow L_{1} &=&
\frac{1}{2}g_{ij}(\f_{0})d\f^{i}_{\perp} * d\f^{j}_{\perp} +
\eps_{IJ}(\f_{0}) \chi^{I} * d \eta^{J}_{\perp} \nonumber \\
& & \;\;\;\;-
\, \Omega_{IJKL}(\f_{0})T^{J}_{\, \overline{M}}\chi^{I}\, \eta^{L}_{0}
\f^{\overline{M}}_{\perp}\, *d \f^{K}_{\perp} \label{act1} \\
{\mathbf L}_{2}\rightarrow L_{2} &=& \frac{1}{2}\left
( \eps_{IJ}(\f_{0})  \chi^{I} d \chi^{J} +
\frac{1}{3} \Omega_{IJKL}(\f_{0}) \chi^{I}\chi^{J}\chi^{K}
\eta^{L}_{0} \right) .\label{act2}
\eea
The notation in these formulae is as follows. Set $\f_{\perp} \in
\Omega^{0} \left( M, \f_{0}^{*}(TX)\right)$, where the $\f_{0}^{i}$ are the
constant maps and the $\f^{i}_{\perp}$ are required to be orthogonal to
the $\f_{0}^{i}$, that is $\int_{M}  *\, \f^{j}_{\perp} = 0$. The
$\eta^{I}$ are also expanded as, $ \eta^{I}= \eta^{I}_{0}+ \eta^{I}_{\perp}$ ,
where the $\eta^{I}_{0}$ are harmonic 0-forms with coefficients in the fibre
$V_{\f_{0}}$ of the ${\mathrm Sp}(n)$ bundle $V \rightarrow X$ and the
$\eta^{I}_{\perp} $ are orthogonal to these $\int_{M}  * \,
\eta^{J}_{\perp}  = 0$.
Though not indicated in the formulae we will, below, also decompose the
$\chi^{I}$ fields in a similar fashion, $\chi^{I} = \chi^{I}_{0} +
\chi^{I}_{\perp} $ where the $\chi^{I}_{0}$ are harmonic 1-forms with
coefficients in the fibre $V_{\f_{0}}$ and the $\chi^{I}_{\perp}$ are
orthogonal to these in the obvious way.

The theory that we will analyze in the following sections is the one
defined in terms of $L_{1}+L_{2}$. This theory is rather simple to get
a handle on, as we will see.

\subsection{Path Integral Properties}
Before proceeding we should mention that we will normalize the bosonic
part of the path
integral measure as done in \cite{RW}. This means that on occasion
certain factors of $2\pi$ will make an appearance and those can always
be traced back to our choice of normalization. Somewhat more involved
is the question of sign of the path integral. Different approaches
to fixing this have been explained in \cite{RW} and \cite{T}, and we
will take the signs to be as given in those references. The question
of framing in the path integral is not adressed here. The
issues involved are spelt out in \cite{RW}.

\subsection{Relationship with Chern-Simons Theory}

In this section we review the relationship between Chern-Simons
theory and the Rozansky-Witten model. Though this relationship has already
been explained in \cite{RW}, we include it here so that we may refer
back to it as we go along.

Recall that the Chern-Simons action is
\be
L_{{\mathrm CS}} = \Tr{ \,\left( AdA \, +\, \frac{2}{3}A^{3}\right)}
\ee
where the trace is understood to be normalized so as to agree with the
standard inner product on the Lie algebra. We compare this with
(\ref{act2}). Notice that there is almost a direct match if we make
the following substitutions
\bea
A^{a} &\rightarrow& \chi^{I} \nonumber \\
\Tr{T_{a}T_{b}} & \rightarrow & \eps_{IJ} \nonumber \\
f_{abc} &\rightarrow& \Omega_{IJKL}(\f_{0})\eta^{L}_{0} ,
\eea
where $T_{a}$ are the generators of the Lie algebra. Also note that
the symmetry properties of the various objects are
reversed. $\Tr {T_{a}T_{b}}$ is symmetric in its arguments while
$\eps_{IJ}$ is antisymmetric, $f_{abc}$ is totally antisymmetric while
$ \Omega_{IJKL}(\f_{0})\eta^{L}_{0}$ is totally symmetric. This is as
it should be since $A^{a}$ is an anti-commuting object while $\chi^{I}$
is commuting. 

In any gauge theory, before performing a perturbative expansion, one needs to
gauge fix, that is to pick a section locally in the space of
connections. In Chern-Simons theory, since we are on
a 3-manifold, we have the trivial connection to use as origin of
the affine space of
connections. A reasonable gauge choice, about the trivial connection, is then
\be
d* A = 0
\ee
which is implemented in the path integral by a delta function
constraint\footnote{Recall that in dimension 1, $\int \ex{ipx} \, dp =
2\pi \d (x) $, is an integral representation of the Dirac delta
`function'. In a lattice approximation of $M$ (\ref{dd}) is understood as
$\prod_{x \in M} \, \d \left( *d*A (x) \right)$ where the product is
over all nodes of the lattice.} 
\be
\d(*d*A) = \int Dt \, \exp{\left( i \int_{M} t\, d*A \right)} \label{dd}
\ee
One should compare this with the second term in (\ref{act1}), with
$t\rightarrow \eta^{I}$. Furthermore (see our review in section
\ref{pirs}), when gauge fixing, in order to balance measures, one must
also introduce the so called Fadeev-Popov ghosts, $c$ and $\bar{c}$,
into the path integral. These, Lie-algebra valued Grassmann odd
zero-forms, enter in
the action as
\be
\int_{M} \Tr{ \left( \bar{c}\, d * (d + A)\, c \right)} .
\ee
Now compare this with (\ref{act1}). The correspondence is readily seen
to be
\bea
c^{a} & \rightarrow & \f^{I} \nonumber \\
\bar{c}^{a} & \rightarrow & \overline{\f}^{\overline{I}} .
\eea

Moreover, we will see in section \ref{pirs} that the topological supersymmetry
of the Rozansky-Witten theory is the natural analogue of the BRST
symmetry of Chern-Simons theory. Even more is true. There is
in any gauge theory an ${\mathrm Sp}(1)$'s worth of BRST symmetry \cite{DJ},
which comes by exchanging the r\^{o}le of the ghosts $c$ and
$\bar{c}$ and this ${\mathrm Sp}(1)$ goes over to the ${\mathrm
Sp}(1)$ of the Rozansky-Witten theory\footnote{One also expects that the
peculiar supersymmetry that exists in Chern-Simons theory \cite{BRT}
also holds here. It would correspond to 
\be
\d \overline{\f}^{\overline{I}} = \eps^{\mu}\, T^{\overline{I}}_{I}\, \chi^{I}_{\mu}
, \;\;\;\; \d \chi^{I}_{\mu} = \varepsilon_{\mu \nu \lambda }\, \eps^{\nu}
\partial^{\lambda} \f^{I} , \;\;\; \d \eta^{I} = \eps^{\mu} D_{\mu} \f^{I}
.
\ee
and a casual glance at the action seems to show that indeed the
symmetry is present, at least in the case of flat space.}. The analogy
is even more remarkable when one notes that the symmetric gauge fixing
is in fact implemented by adding, 
\be
Q \overline{Q} \, \left( \Tr \,  A * A \right),
\ee
to the action, which should be compared with (\ref{dbrst}).

Given the intimate relationship with Chern-Simons one would expect
that both the IHX and the AS relations would hold in the
Rozansky-Witten theory. That the IHX relation is satisfied was
established in \cite{RW} and corresponds to the geometric analogue of
the Jacobi identity, namely to the Bianchi identity for the $\Omega$
tensor. At a naive level the AS relation does not appear to be true in
Rozansky-Witten theory, however it holds in the most meaningful
way.  

\begin{figure}[h]
\centerline{\epsffile{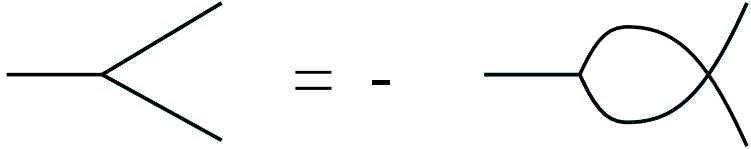}}
\caption{The AS Relation}\label{asrel}
\end{figure}
\begin{figure}[h]
\centerline{\epsffile{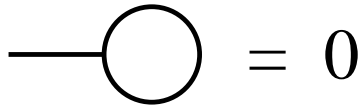}}
\caption{The Tadpole Diagram}\label{tad}
\end{figure}

Recall that the AS relation in Chern-Simons theory amounts to the
anti-symmetry property of the structure constants, $f_{abc}$, of the
Lie algebra, that is, $f_{abc} = -  f_{acb}$, as depicted in
Figure~\ref{asrel}. On the other hand, in the Rozansky-Witten theory
the ``structure
constant'' which appears in diagrams is the completely symmetric
tensor $\Omega_{IJKL}$ and so the identity implied by
Figure~\ref{asrel} appears
to be violated. However, one must recall that, in reality, the
vertices in the Rozansky-Witten theory are connected to Grassmann odd
objects and that if one thinks of the vertices as incorporating this
Grassmann character then the vertex obeys the AS relation. So for
example in the Rozansky-Witten theory, any diagram which has a loop
centered at a vertex, as shown in Figure~\ref{tad}, vanishes which is a fact
completely consistent with the AS relation. The reason it vanishes is
that while $\Omega$ is symmetric, the labels in the loop are
contracted by $\eps^{IJ}$ which is anti-symmetric and that the
contraction is anti-symmetric is due to the fact that we had Grassmann
odd variables there. 
\noindent In section \ref{as} we
will be quite explicit about how the AS relation arises for $b_{1}(M)=3$.

For all the similarity there is one important difference between the
two theories. The vertex in the Rozansky-Witten model carries a
Grassmann odd harmonic mode, $\eta^{I}_{0}$. This means that this
vertex may never appear more than $2n$ times in any diagram. Thus 
there is a cut-off built into the perturbative expansion of the
Rozansky-Witten theory.

\subsection{Compact and non-Compact $X$}\label{X}

The spaces $X_{G}$ that are associated to a group $G$ really
correspond to certain moduli spaces of monopoles. These spaces are
hyper-K\"{a}hler but non-compact. Nevertheless, they are asymptotically
flat. The dependence that one finds on $X_{G}$ is through terms of the
form
\be
\int_{X_{G}} \Tr \left( \frac{{\mathrm
R(X_{G})}}{2\pi}\right)^{2j_{1}}\dots \, \Tr \left( \frac{{\mathrm
R(X_{G})}}{2\pi}\right)^{2j_{n}}, \label{inttypes}
\ee
with $\sum_{i} j_{i} =n$, or with explicit dependence on the
holomorphic 2-form, such that the
integrand is a top form. Since the manifolds are
asymptotically flat integrals of this type will make sense. To be sure
that the invariants do not trivially vanish we need to know if
integrals of the form (\ref{inttypes}) are zero or not.

Non-compact hyper-K\"{a}hler manifolds abound. Examples include the
Atiyah-Hitchin manifold $X_{AH}$, which is the $SU(2)$ 2-monopole
moduli space as well as $T^{*}{\Bbb CP}^{n}$ for which Calabi \cite{C}
exhibited hyper-K\"{a}hler
metrics. More
generally, one
has a procedure for producing examples. Suppose that
one is given a hyper-K\"{a}hler manifold $X$ admitting a Lie group action of
isometries $G$ which preserve the hyper-K\"{a}hler structure with $\mu : X
\rightarrow \lg^{*} \times {\Bbb R}^{3}$ the corresponding moment
map, where $\lg^{*}$ is the dual of the Lie algebra of $G$. Then, the
hyper-K\"{a}hler quotient construction \cite{HKLR} guarantees that if
$G$ acts freely on $\mu^{-1}(0)$ with a Hausdorff quotient then the
quotient manifold $\mu^{-1}(0)/G$ (denoted $X//G$) is once more
hyper-K\"{a}hler (with
the hyper-K\"{a}hler metric being the induced one). Starting with the
hyper-K\"{a}hler manifold ${\Bbb C}^{m} \times {\Bbb C}^{m}$ and
quotienting with various groups gives rise to many known examples of
hyper-K\"{a}hler manifolds including the monopole moduli spaces of
interest. For $n>1$ we are unaware of any calculations for integrals
of the form (\ref{inttypes}). There appears to be a dearth of
information on the properties of such integrals. We will proceed under
the assumption that there are sufficiently many manifolds for which
integrals of the type (\ref{inttypes}) are finite and non-vanishing.

Of course, once one has the invariants at one's disposal, they are
defined for any hyper-K\"{a}hler $X$ and not just $X_{G}$. So that, in
particular, one may also consider compact manifolds. However, while
non-compact hyper-K\"{a}hler manifolds are plentiful the compact variety are
rare birds indeed. There are essentially two series of examples
\cite{Be}. The first is
made up of a resolution of the n-fold symmetric product of $K3$ surfaces and
is denoted by $S^{[n]}$. The Douady space, $S^{[n]}$, is a (real) $4n$
dimensional irreducible hyper-K\"{a}hler manifold. The second series,
denoted by $K_{n}$, is related to the Douady space,
$T^{[n]}$, of the n-fold symmetric product of the four dimensional
torus $T$.
$T^{[n]}$ is not irreducible while $K_{n}$ is. There is only one known
example which is neither of type $S^{[n]}$ or $K_{n}$.

What we would really like is to get a handle on integrals of
the form (\ref{inttypes}). Fortunately, very recently, computations of
the even Chern numbers for the $S^{[n]}$ series have been made for
$n=1, \dots
, 7$.\footnote{We thank L. G\"{o}ttsche for making these computations
available to us.} It is quite remarkable that these are all non-zero and
positive. One can check to see that the Chern characters of the
tangent bundle in these cases do not vanish. As far as we are aware
there are  no similar computations
available for the $K_{n}$ series except for the Euler characteristic,
which is again strictly positive and again due to  L. G\"{o}ttsche
\cite{G}.

\subsection{Product Groups (Manifolds)}

Since the Generalized Casson invariant can be morally viewed as the Euler
characteristic of the moduli space of flat $G$ connections, one has
immediately that the invariant for a product group is the product of
the invariant of each group factor. Let the group have the form $G=
G_{1} \times G_{2}$, then ${\mathrm Hom}(\pi_{1}(M), G)/G= {\mathrm
Hom}(\pi_{1}(M), G_{1})/G_{1} \times {\mathrm Hom}(\pi_{1}(M), G_{2})/
G_{2}$, or ${\mathcal M}(G)= {\mathcal M}(G_{1})\times {\mathcal M}(G_{2})$
and consequently $\chi ({\mathcal M}(G))= \chi({\mathcal M}(G_{1}))\,
.\,  \chi({\mathcal M}(G_{2}))$.

How does the Rozansky-Witten invariant behave when we consider product
groups? The answer is that it factorizes as it should. To pass from
the gauge theory to the sigma model one uses the dictionary $G
\rightarrow X_{G}$, which for products reads $G_{1}\times G_{2}
\rightarrow X_{G_{1}} \times X_{G_{2}}$. The hyper-K\"{a}hler
structure of a product manifold is the natural product
hyper-K\"{a}hler structure. The path integral factorizes since the
space of maps factorizes and the Lagrangians split into the sum of two
pieces, one of which only involves objects associated with $X_{G_{1}}$, the
other involving objects only depending on $X_{G_{2}}$. We have then
that
\be
Z_{X_{1} \times X_{2}}^{RW}[M] = Z_{X_{1}}^{RW}[M] \, . \,
Z_{X_{2}}^{RW}[M] . \label{factor}
\ee

This simple observation has immediate consequences for the invariant,
if the dependence on $X$ is only
through characteristic classes. For, if this is the case, then we may
expand the partition function as
\be
Z_{X}^{RW}[M] = \int_{X}\sum_{j_{1}\geq \dots \geq j_{n}}^{n} \,
\Tr \left( \frac{{\mathrm
R(X)}}{2\pi}\right)^{2j_{1}}\dots \, \Tr \left( \frac{{\mathrm
R(X)}}{2\pi}\right)^{2j_{n}} \, I_{(j_{1}, \dots , j_{n})}^{n}(M)
. \label{pfunct}
\ee
(It is understood that in (\ref{pfunct}) one picks out the form of
degree $4n$ in the integrand, that is $\sum_{k=1}^{n}j_{k}=
n$. The actual form of the integrals $I^{n}_{j_{1}, \dots , j_{n}}(M)$
depends very much on the first Betti number of $M$.) Now the
factorization property
(\ref{factor}) implies that, in fact,
\be
Z_{X}^{RW}[M] = \int_{X} \, \ex{ \sum_{j=1}^{n} \Tr \left
( \frac{{\mathrm R}(X)}{2\pi} \right)^{2j}\, I_{j}(M) } , \label{expf}
\ee
where
\be
I_{j}(M) = I^{j}_{(j, 0 ,\dots ,0)}(M) .
\ee

Consequently, the $I_{(j_{1}, \dots , j_{n})}^{n}(M)$ are completely
determined by the $I_{j}(M)$ for $ j\leq n$.\footnote{One might have
thought that by suitably juggling terms proportional
to, 
\be
\Tr \left
 ( \frac{{\mathrm R}(X)}{2\pi} \right)^{2j} \, \Tr \left
( \frac{{\mathrm R}(X)}{2\pi} \right)^{2p}
\ee
with $j+p \leq n$ in the
exponent, one might still be able to satisfy (\ref{factor}). However,
since we can find hyper-K\"{a}hler manifolds $X_{j}$, $X_{p}$ and
$X_{4n-4j-4p}$ of dimension $4j$, $4p$ and $4(n-j-p)$ respectively,
such a term would spoil the factorization property. In the text we
will see explicitly that the partition function indeed
takes the form of (\ref{expf}), when $b_{1}(M) \geq1$.} This has quite
drastic implications. Since at any $n$ we are claiming
that there is only one new integral $I_{n}(M)$ that arises, there is
then only at most one new invariant at the given $n$. Consequently there are
only, maximally, a ${\Bbb Z}$'s worth of invariants! We will see, in
the following sections, that if $M$ has rank $\geq 1$ then both the
hypothesis that $X$ enters only through its Chern numbers and the
conclusions drawn hold.

The interesting case then is $b_{1}(M)=0$. In this case we cannot show
that the $X$
only enters through its Chern numbers. This is just as well since it
is believed that the number of LMO invariants grows rather more rapidly than
linearly with dimension (degree). However, for $n=1$, Rozansky and
Witten showed that the invariant is proportional to $\Tr R^{2}$ and
one can also show that the double theta at $n=2$ 
\begin{figure}[h]
\centerline{\epsffile{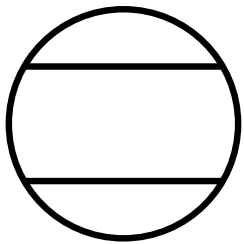}}
\caption{The Double Theta Diagram}\label{dtheta}
\end{figure}
\noindent is proportional to $\Tr R^{4}$. Since the Mercedes Benz diagram 
\begin{figure}[h]
\centerline{\epsffile{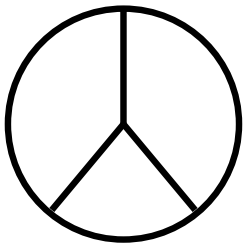}}
\caption{The Mercedes Benz Diagram}
\end{figure}
\noindent is proportional to the
double theta this means that also at $n=2$ there is only one new invariant.

While the main thrust of our physical computations is to avoid working
directly with diagrams, one aspect of the factorization property for
any $M$ is very simple to describe in terms of the diagrammatic
expansion. One deduces from (\ref{factor}) that for product manifolds the
connected diagrams vanish while the product diagrams factor to
reproduce the formula. The vanishing of the connected diagrams is a
simple consequence that one gets from considering how the 
$\eta^{I}$ zero modes enter into the diagram. 

On a product manifold,
$X = X_{1}\times X_{2}$, with ${\mathrm dim}_{{\Bbb R}}X_{1} = 4n_{1}$,
${\mathrm dim}_{{\Bbb R}}X_{2} = 4n_{2}$ and $n_{1}+n_{2} = n$,
the 2-form $\eps_{IJ}$ becomes the sum of the
holomorphic symplectic 2-forms of each factor, i.e. $\eps_{IJ}=
\eps_{IJ}^{1} +  \eps_{IJ}^{2}$. Likewise, the curvature tensor splits
as $\Omega_{IJKL}= \Omega_{IJKL}^{1}+ \Omega_{IJKL}^{2}$. However, one
must remember that in diagrams the curvatures are connected by
$\eps^{IJ}$ (coming from propagators). Which means that vertices with
$\Omega_{IJKL}^{1}$ assigned to them can only be connected to other vertices
with $\Omega_{IJKL}^{1}$ assigned to them since the $\eps^{IJ}$ do not
`mix' manifolds. This, in turn, means that for connected diagrams
one has assigned $\Omega_{IJKL}^{1}$ to every vertex or one has
assigned $\Omega_{IJKL}^{2}$ to every vertex. There are $2n_{1}$
harmonic $\eta^{I}$ modes, denoted $\eta^{I}_{1}$, from $X_{1}$ and $2n_{2}$
harmonic $\eta^{I}$ modes, denoted $\eta^{I}_{2}$, from $X_{2}$.

Hence, in any given
connected diagram with $m$ vertices one has exactly $m$ $\eta^{I}_{1}$
or $m$ $\eta^{I}_{2}$ zero modes appearing. In order to get a
non-vanishing answer for the integral over the harmonic modes the
product of connected diagrams appearing in one Feynman diagram must be
such that exactly $2n_{1}$ vertices can have $\Omega_{IJKL}^{1}$
assigned to them and $2n_{2}$ vertices can have $\Omega_{IJKL}^{2}$
assigned to them. For a completely
connected Feynman diagram, with $2n$ vertices, this is not possible if
both $n_{1}$ and $n_{2}$ are non-zero.

ÿÿ

\subsection{Observables and $Q$-cohomology}

There are a number of observables that can be defined. The expectation
value of each of these is potentially a new invariant, though, as we
will see, they may be invariants that we have already encountered.

For us, the basic set of observables involves powers of the
holomorphic symplectic 2-form $\eps$, taken at some point of $M$ by
pull-back. However as will be seen below, the precise point on $M$ at
which the form is evaluated is
immaterial and therefore will be suppressed from the
notation. Let\footnote{That the following observables make sense can
be seen by noting that they should be viewed as the pull back of the
evaluation of the k-th wedge product of the holomorphic 2-form on $2k$
(Grassmann odd) tangent vectors.}
\be
{\mathcal O}(k) = d_{k}
\left(-\frac{1}{2}\eta^{I}\eps_{IJ}\eta^{J}\right)^{k} .
\ee

The coefficients $d_{k}$ will be fixed below. The invariants
of the manifold $M$ are defined by a path integral which has an
insertion of ${\mathcal O}(k)$, that is
\be
Z^{RW}_{X}[M, {\mathcal O}(k)] = \int D\f D\eta D\chi \, \ex{-S} \,
{\mathcal O}(k) .
\ee

Let $< A >$ denotes the expectation value of $A$, with respect to some
set of fields $\Phi$ and some action $S(\Phi)$,
\be
< A > = \int D\Phi\, \ex{-S(\Phi)} \, A .
\ee
The measure on the harmonic $\eta^{I}$ modes is determined by
\be
<\eta^{I_{1}} \, \dots \, \eta^{I_{2n}}> = \eps^{I_{1}\, \dots \, I_{2n}} ,
\ee
where we mean that one integrates only over the harmonic $\eta^{I}$
modes and for which the action is taken to be zero.

We, partially, fix the coefficients $d_{k}$ by demanding that
\be
d_{0}=1,
\ee
and that
\be
\frac{1}{(2\pi)^{2n}} \int_{X}\sqrt{g} \, d^{4n} \f_{0} \, <{\mathcal
O}(k) {\mathcal O}(n-k)>= 1  , \label{norm}
\ee
where $\sqrt{g} \, d^{4n} \f_{0}$ is the Riemannian measure on
$X$. Notice that when ${\mathrm dim}_{{\Bbb R}}X = 4$, that this
specifies the value
of $d_{1}$, while for ${\mathrm dim}_{{\Bbb R}}X = 8$, $d_{0}$, and $d_{2}$ are
determined and $d_{1}$ is fixed up to a sign.

One of the most important properties of this class of observables is
that
\be
\int D\f D\eta D\chi \ex{-S} \; {\mathcal O}(n) = \left(|H_{1}(M, {\Bbb
Z})| \right)^{n}.\label{normn}
\ee
This follows from the normalization that we have chosen in
(\ref{norm}) as well as the observation in \cite{RW} that the
Ray-Singer torsion provides a natural volume form which includes the
Riemannian volume of $X$ times $|H_{1}(M, {\Bbb Z})|' = |{\mathrm Tor}
H_{1}(M, {\Bbb Z})|$. The reason
that it is $|H_{1}(M, {\Bbb Z})|$ rather than $|H_{1}(M, {\Bbb Z})|'$
that appears in (\ref{normn}) is that the observables, ${\mathcal
O}(k)$, vanish for manifolds that are not $\QHS$'s (this follows by a
count of vertices similar to those made in section \ref{sec.proof2}). 

We can take the coefficients $d_{k}$, to be
\be
d_{k} = \frac{(2\pi)^{2k}}{\left({\mathrm Vol}(X)n!\right)^{k/n}} =
\left(d_{n}\right)^{k/n} ,
\ee
so that the observables enjoy the following property
\be
{\mathcal O}(p){\mathcal O}(q) = {\mathcal O}(p+q).
\ee

We should explain why these are good operators to consider in the
theory. Let ${\cal O}(\omega)$ be defined by
\be
{\cal O}(\omega) = \omega_{I_{1} \dots I_{k}} \, \eta^{I_{1}} \dots
\eta^{I_{k}} 
\ee
where the $\omega_{I_{1} \dots I_{k}}$ are both the components of a
$\partial$ closed $(k,0)$ form as well as being the components of a
$\overline
{\partial}$ closed $(0,k)$ form. (Here $\partial$ and
$\overline{\partial}$ are the holomorphic and anti-holomorphic
Dolbeault operators on $X$ and the correspondence is given by the
isomorphism between $TX^{(1,0)}$ and $TX^{(0,1)}$.) In equations this means
that we want
\be
\partial_{I} \, \omega_{I_{1} \dots I_{k}} dz^{I}dz^{I_{1}} \dots
dz^{I_{k}} = 0 ,
\ee
and
\be
\overline{\partial}_{\overline{I}}\,  \omega_{\overline{I}_{1} \dots
\overline{I}_{k}} d\bar{z}^{\overline{I}}d\bar{z}^{\overline{I}_{1}} \dots
d\bar{z}^{\overline{I}_{k}} = 0 ,
\ee
where
\be
\omega_{\overline{I}_{1} \dots \overline{I}_{k}} =
\omega_{I_{1}  \dots I_{k}} \, T^{I_{1}}_{\; \overline{I}_{1}} \dots
T^{I_{k}}_{\; \overline{I}_{k}},
\ee
and $T^{I}_{\; \overline{I}} = g_{\overline{I} J}\eps^{JI}$ (see the
appendix for more details about hyper-K\"{a}hler manifolds and for our
conventions).

The important property of such operators is that they are both $Q$ and
$\overline{Q}$ closed. That is,
\bea
Q \, {\cal O}(\omega) & = & Q \f^{i} \, \frac{\partial {\cal
O}(\omega)}{\partial \f^{i}} \nonumber \\
& = & \eta^{I} \, \partial_{I} \omega_{I_{1} \dots I_{k}} \,
\eta^{I_{1}} \dots \eta^{I_{k}} \nonumber \\
& = & 0 ,
\eea
where the third equality follows from the fact that $\omega$ is
$\partial$ closed.
(While these operators are $Q$ closed they are not $Q$ exact, if
$\omega$ is non-trivial in cohomology, since $Q$ exactness would imply
$\partial $ exactness of $\omega$.) Similarly
\bea
\overline{Q} \, {\mathcal O}(\omega) & = & T^{\bar{J}}_{\; J}\eta^{J} \,
\partial_{\bar{J}}{\mathcal O}(\omega) \nonumber \\
& = & \eta^{I} \, \partial_{I} \omega_{I_{1} \dots I_{k}} \,
\eta^{I_{1}} \dots \eta^{I_{k}} \nonumber \\
& = & 0 ,
\eea

Another important property of such observables is that they are
essentially $d$ closed as well, where $d$ is the exterior derivative
on $M$. Here essentially means that this holds because the path
integral is concentrated along the constant maps. But
since the dependence of the observables on $M$ is via pullback with
respect to $\f$, they do not depend on the point at which they sit on
$M$.

\section{Outline of the Proof of Heuristic Theorem 2.} \label{sec.proof2}
The strategy of the proof will be to decide which types of
Feynman diagrams can contribute and then to find a way of encoding all the
relevant information without doing any expansions. The one piece of
information that we will use continuously is that, on expanding out
the interaction part of the action, the interaction terms will be of
the form
\be
V_{1}^{p}\, V_{2}^{q} , \label{vpower}
\ee
where
\bea
V_{1} & = &\frac{1}{6} \int_{M}  \Omega_{IJKL}(\f_{0})
\chi^{I}\chi^{J} \chi^{K} \eta^{L}_{0} \\
V_{2} &= &- \int_{M}
\Omega_{IJKL}(\f_{0})T^{J}_{\, \overline{M}}\chi^{I}\, \eta^{L}_{0}
\f^{\overline{M}}_{\perp}\, *d \f^{K}_{\perp} .
\eea
There is a condition on the number of interactions that
arise from the fact that both vertices $V_{1}$ and $V_{2}$ are
linear in $\eta^{I}_{0}$, which we recall is constant on $M$. Since
$\eta^{I}_{0}$ is constant on $M$, this means that the ``path integral''
over this ``field'' is actually a 2n-fold product of Berezin
integrals (the exact specification of the measure is described in
appendix 
\ref{sec.manipulations}). Furthermore, from the rules that we describe
for such variables,
we see that the integration will vanish identically unless the
integrand includes the product of all $2n$ components of
$\eta^{I}_{0}$ (each occuring exactly once). Consequently one has
\be
p+q=2n . \label{p+q}
\ee

The rest of the proof depends on how many $\chi^{I}$ harmonic modes
there are. The importance of these modes lies in the fact that, like
the $\eta^{I}_{0}$, they will only appear in the
vertices.\footnote{The harmonic modes do not appear in the quadratic
terms, for example $\int_{M} \eps_{IJ}(\f_{0})\chi^{I}* d \eta^{J}_{\perp}
=\int_{M} \eps_{IJ}(\f_{0})\chi^{I}_{\perp}* d \eta^{J}_{\perp}$, by an
integration by parts.} Hence, for the
same reason as for the harmonic $\eta^{I}$, the path integral will vanish
identically unless the integrand includes the product of all the
harmonic components of $\chi^{I}$. The number of harmonic modes of
$\chi^{I}$ is $2n b_{1}(M)$. This means that there is another condition
that must be satisfied to ensure that the integral has a chance of not
vanishing, which is
\be
3p + q \geq 2n b_{1}(M). \label{p>b}
\ee
This inequality comes by noticing that if $3p+q <2n b_{1}(M)$ then
certainly the integrand will not have the required product of
harmonic $\chi^{I}$. As the constraint (\ref{p>b})
depends on $b_{1}(M)$, we go through the
possible $b_{1}(M)$ values and along the way we will strengthen it.
Subtracting (\ref{p+q}) from (\ref{p>b}) gives the constraint
\be
p \geq n (b_{1}(M)-1) . \label{pb}
\ee

\subsection{Manifolds with $b_{1}(M)\geq 4$}

From (\ref{pb}) we see that if $b_{1}(M)\geq 4$, then clearly $p > 2n$, but
then (\ref{p+q}) certainly cannot be satisfied. Consequently, we see
that it is in fact impossible to integrate up the harmonic modes of
$\chi^{I}$ and hence the Rozansky-Witten path integral vanishes.

\section{Manifolds with $b_{1}(M)=3$}

The inequality (\ref{pb}) says, for $b_{1}(M)=3$, that $p\geq 2n$ and
this intersects with (\ref{p+q}) only if $p=2n$ and $q=0$. This
condition tells us that we are to ignore $V_{2}$ completely, so that
perturbatively one is interested in
\be
V_{1}^{2n} . \label{intb3}
\ee
However, we know more. Since there are $6n$ $\chi^{I}$ harmonic modes,
all the $\chi^{I}$ appearing in (\ref{intb3}) must be harmonic. This
means that vertex $V_{1}$ effectively reduces to
\be
V_{1} = \frac{1}{6}\int_{M} \,  \Omega_{IJKL}(\f_{0}) \chi^{I}_{0}
\chi^{J}_{0}\chi^{K}_{0}\eta^{L}_{0} \label{vert}
\ee
and in turn that the Lagrangians (\ref{act1}, \ref{act2}) reduces
(with a slight rearrangement) to
\bea
L_{\perp} &=&
\frac{1}{2}g_{ij}(\f_{0})d\f^{i}_{\perp}*d\f^{j}_{\perp} +
\eps_{IJ}(\f_{0}) \chi^{I}_{\perp}*d \eta^{J}_{\perp} +
\frac{1}{2}\eps_{IJ}( \f_{0}) \chi^{I}_{\perp} d
\chi^{J}_{\perp} \label{nact1} \\
L_{0} &=&
\frac{1}{6} \Omega_{IJKL}(\f_{0}) \chi^{I}_{0}\chi^{J}_{0}\chi^{K}_{0}
\eta^{L}_{0}  \label{nact2}
\eea
where a zero subscript indicates the field is harmonic while the
$\perp$ subscript indicates that the field is orthogonal to the
harmonic modes on $M$. The path integral to be performed is,
symbolically,
\be
\int d\Phi_{0} \; \ex{ -\int_{M} L_{0} } \;  \int D\Phi_{\perp}
\; \ex{-\int_{M} \, L_{\perp}} ,
\ee
where $\Phi$ denotes the set of fields. Rozansky and Witten \cite{RW}
have shown that
\be
\int D\Phi_{\perp} \; \ex{-\int_{M} \, L_{\perp}} = \left(|{\mathrm
H}_{1}(M, {\Bbb Z})|' \right)^{n},
\ee
where $|{\mathrm H}_{1}(M, {\Bbb Z})|' $ is the order of ${\mathrm Tor
\, H}_{1}(M,{\Bbb Z})$. Consequently, the Rozansky-Witten invariant in
this case has the very succinct representation
\be
Z^{RW}_{X}[M]= \left(|{\mathrm
H}_{1}(M, {\Bbb Z})|' \right)^{n}\int d\Phi_{0} \; \ex{ -\int_{M}
L_{0} }. \label{lefttodo}
\ee
Let $\o^{\a}$, $\a = 1, \, 2, \, 3$ be a basis of $H^{1}(M, {\Bbb
Z})$. Since $\chi_{0}^{I}$ is harmonic on $M$, we see that it must have
the expansion $\chi_{0}^{I} = \chi^{I}_{\a} \o^{\a}$, where the
coefficients $\chi^{I}_{\a}$, for each $\alpha = 1,2,3$, are
generators of $\wedge \f^{*}T^{(1, 0)}X$. Substitution of this
expansion in (\ref{vert}) gives
\bea
V_{1} &= &  \frac{1}{6} \Omega_{IJKL}(\f_{0}) \chi^{I}_{\a}
\chi^{J}_{\b}\chi^{K}_{\gg}\eta^{L}_{0} \int_{M} \,
\o^{\a}\o^{\b}\o^{\gg} \nonumber \\
&= & \frac{1}{6} \Omega_{IJKL}(\f_{0}) \chi^{I}_{\a}
\chi^{J}_{\b}\chi^{K}_{\gg}\eta^{L}_{0} \eps^{\a \b \gg}
\; I(M) , \label{vert1}
\eea
where
\be
I(M) = \frac{1}{6}\eps_{\a \b \gg} \int_{M} \,
\o^{\a}\o^{\b}\o^{\gg} .
\ee
By the arguments that we presented at the start of this section
culminating in (\ref{intb3}) we immediately have that the
Rozansky-Witten invariant is proportional to $I(M)^{2n}$.

To determine the coefficient, it suffices to compute the
invariant for any 3-manifold of rank 3
as it does not depend on $M$ but only on $X$. In \cite{T} (equation (3.23)),
the second author showed that the
Rozansky-Witten invariant for the 3-torus, $T^{3}$, is equal to the Euler
characteristic of $X$ if it is compact. (More generally, it is equal to the
integral of the Euler form over $X$) Denote this (in
both cases) by ${\mathbf e}(X)$. We also have $I(T^{3})^{2} = 1$
so that we find for manifolds
$M$, with $b_{1}(M)=3$,
\be
Z^{RW}_{X}[M] = {\mathbf e}(X)\left(|H_{1}(M, {\Bbb Z})|'  I(M)^{2}
\right)^{n} .
\ee

One may relate this back to the Lescop invariant $\lambda (M)$, as
\be
I(M)^{2} = \frac{\lambda (M)}{|H_{1}(M, {\Bbb Z})|' }
\ee
so that
\be
Z^{RW}_{X}[M] = c_{X} \lambda (M)^{n} , \label{b3}
\ee
where
\be
c_{X} = {\mathbf e}(X) .
\ee

One can also get this result by performing the finite dimensional
integrals that were left to be done in (\ref{lefttodo}). Indeed these
are computed directly from (\ref{c2}).

As already explained in section \ref{X} there are compact and
non-compact $X$
for which ${\mathbf e}(X)$ is non-vanishing for every
$n$, so that (\ref{b3}) is not empty. For example, for the $S^{[n]}$
series one has the generating function
\be
\sum_{n=0}^{\infty} {\mathbf e}(S^{[n]}) t^{n} = \prod_{k=1}^{\infty}
\frac{1}{\left( 1 - t^{k} \right)^{24}} , \label{seul}
\ee
while for the $K_{n}$ series one has
\be
{\mathbf e}(K_{n}) = (n+1)^{3} \sigma_{1}(n+1) ,
\ee 
where $\sigma_{1}(n)$ is the sum of the divisors of $n$. 
It is amusing that, for $S^{[n]}$, if one replaces the Casson
invariant with the indeterminant $t$, then summing over $n$ in
(\ref{b3}) reproduces the
generating function for the Euler characteristics (\ref{seul}). 

On the non-compact side $\left(X_{AH}\right)^{n}$ will do, since
${\mathbf e}(X_{AH}^{n}) = {\mathbf e}(X_{AH})^{n}=1$.

\subsection{The AS Relation}\label{as}
The vertex with the $\chi$ zero modes attached is
\be
V_{1} =   \frac{1}{6} \Omega_{IJKL}(\f_{0}) \chi^{I}_{\a}
\chi^{J}_{\b}\chi^{K}_{\gg}\eta^{L}_{0} \int_{M} \,
\o^{\a}\o^{\b}\o^{\gg} .
\ee
If one extracts the part that depends on $X$ from the dependence on
$M$ we can write the vertex as
\be
V_{1} = F_{\a \b \gg}(X) \, \int_{M} \,
\o^{\a}\o^{\b}\o^{\gg} ,
\ee
where $F_{\a \b \gg}(X) = \frac{1}{6} \Omega_{IJKL}(\f_{0}) \chi^{I}_{\a}
\chi^{J}_{\b}\chi^{K}_{\gg}\eta^{L}_{0}$, is totally antisymmetric
in its three labels. Consequently, it is the vertex $F_{\a \b \gg}(X)$
that satisfies the AS relation. In Chern-Simons theory a similar
vertex arises when one replaces
the gauge connection with harmonic modes in the cubic term, that is,
one sets $A^{a} = A^{a}_{\a} \o^{\a}$. The cubic term is now
proportional to $K_{\a \b \gg}(G)\, \int_{M}
\,\o^{\a}\o^{\b}\o^{\gg}$, where $K_{\a \b \gg}(G) = f_{abc}
A^{a}_{\a} A^{b}_{\b} A^{c}_{\gg} $. Notice that $K_{\a \b \gg}(G)$ is
also antisymmetric in its three labels. The antisymmetry of $F_{\a \b
\gg}(X)$ is a consequence of the symmetry of $\Omega$ and the
anticommuting properties of the Grassmann variables $\chi^{I}_{\a}$. The
antisymmetry of $K_{\a \b \gg}(G)$ rests on the antisymmetry of the
structure constants of $G$, $f_{abc}$ and the fact that the
variables $A^{a}_{\a}$ commute with each other.

\section{Manifolds with $b_{1}(M)=2$}

In this case from (\ref{pb}) we learn that $p \geq n$. However, we can
show that $q=0$. To see this, note that in the vertex $V_{1}$
there can be at most two harmonic $\chi^{I}$, since the wedge product of
three would be zero ($b_{1}(M)=2$). This means that we can refine
(\ref{p>b}) to obtain the inequality
\be
2p + q \geq 4n
\ee
which together on subtracting (\ref{p+q}) tells us that $p \geq
2n$. Hence, once more we find that $p=2n$ and $q=0$.
What this means for us is that we may ignore $V_{2}$ and also in order
to guarantee that the harmonic modes are accounted for, two and only
two of the $\chi^{I}$ appearing in $V_{1}$ must be harmonic. One sets
\be
V_{1} = \frac{1}{2}\int_{M} \,  \Omega_{IJKL}(\f_{0}) \chi^{I}_{\perp}
\chi^{J}_{0}\chi^{K}_{0}\eta^{L}_{0} .\label{vertb2}
\ee
Let $\o^{\a}$, $\a = 1, \, 2$, be a basis of $H^{1}(M, {\Bbb
Z})$. Since $\chi_{0}^{I}$ is harmonic on $M$, we see that it must have
the expansion $\chi_{0}^{I} = \chi^{I}_{\a} \o^{\a}$,
where the coefficients $\chi^{I}_{\a}$, for each $\alpha=1,2$, are
generators of $\wedge
\f^{*}T^{(1, 0)}X $. Inserting this into (\ref{vertb2}),
yields
\be
V_{1} =-\frac{1}{2}\Omega_{IJKL}(\f_{0}) \chi^{J}_{\a} \chi^{K}_{\b}
\eta^{L}_{0}\int_{M} \, \chi^{I}_{\perp} \o^{\a} \o^{\b} .
\ee
While the wedge product $\o^{\a} \wedge \o^{\b}$ is exact it is not
harmonic. Set $\o^{\a} \wedge \o^{\b} = \eps^{\a \b} \, d g $, where
$g$ is a one form on $M$. $g$ is only defined up to exact pieces
so to be definite we demand that $d*g =0$.

The actions become,
\bea
\int_{M}L_{1} &=&
\int_{M} \, \left(\frac{1}{2}g_{ij}(\f_{0})d\f^{i}_{\perp}*d\f^{j}_{\perp} +
\eps_{IJ}(\f_{0}) \chi^{I}_{\perp}*d \eta^{J}_{\perp}  \right)\label{act1b2} \\
\int_{M}L_{2} &=& \frac{1}{2}\eps_{IJ}( \f_{0}) \int_{M} \chi^{I}_{\perp} d
\chi^{J}_{\perp}
- \Lambda_{I} \int_{M} \, d \chi^{I}_{\perp} \, g . 
\label{act2b2}
\eea
where we have set $
\Lambda_{I} = \frac{1}{2}\Omega_{IJKL}(\f_{0}) \eps^{\a \b}
\chi^{J}_{\a} \chi^{K}_{\b} \eta^{L}_{0}$.

We now complete the square
\be
\int_{M}L_{2} = \frac{1}{2}\eps_{IJ}( \f_{0}) \int_{M}
\left(\chi^{I}_{\perp} +\eps^{IK}\Lambda_{K}g  \right)d
\left(\chi^{J}_{\perp} + \eps^{JL} \Lambda_{L}g \right)
+\frac{1}{2} \Lambda_{I} \eps^{IJ}\Lambda_{J}\int_{M} g dg .
\label{csq}
\ee
Notice also that the part of (\ref{act1b2}) that involves
$\chi_{\perp}^{I}$ is
\be
\int_{M} \eps_{IJ}(\f_{0}) \chi^{I}_{\perp}*d \eta^{J}_{\perp} =
\int_{M} \eps_{IJ}(\f_{0}) \left( \chi^{I}_{\perp} +
\eps^{IK}\Lambda_{K}g  \right) *d \eta^{J}_{\perp} ,
\ee
since we have chosen $d*g=0$. We now change variables in the path
integral. Since $\chi^{I}_{\perp}$ only appears in the action in the
combination $\hat{\chi}^{I}_{\perp}= \chi^{I}_{\perp}
+ \eps^{IK}\Lambda_{K}\,g $, we change variables to
$\hat{\chi}^{I}_{\perp}$. Does such a change of variables make sense?
The answer is yes. Firstly, since $\Lambda_{I}$ is Grassmann odd we are
maintaining the grading of the fields. Secondly, $g$ lives in
$H^{1}_{\perp}(M)$ so this character of the field is also
preserved. Lastly, the Jacobian for such a change of variables
is unity since the object by which we are shifting $\chi^{I}_{\perp}$ does
not depend on $\chi^{I}_{\perp}$.

Hence the actions (\ref{act1b2}) and (\ref{act2b2}) can be grouped as
follows
\be
S_{0} + S_{\perp} = \int_{M} \left( L_{1} + L_{2} \right)
\ee
where
\bea
S_{0} &=&  \frac{1}{2} \Lambda_{I} \eps^{IJ}\Lambda_{J}\int_{M} g dg , \\
S_{\perp} & = & \int_{M} \,
\left(\frac{1}{2}g_{ij}(\f_{0})d\f^{i}_{\perp}*d\f^{j} _{\perp} \right.
\nonumber \\
& & \;\;\; + \left.
\frac{1}{2} \eps_{IJ}( \f_{0}) \int_{M}\hat{\chi}^{I}_{\perp} d
\hat{\chi}^{J}_{\perp} +
\eps_{IJ}(\f_{0}) \hat{\chi}^{I}_{\perp}*d \eta^{J}_{\perp}  \right) .
\eea
This shows us, once more, that the path integral over all the fields
splits nicely as,
\be
\int d\Phi_{0} \; \ex{ -S_{0} } \;  \int D\Phi_{\perp}
\; \ex{-S_{\perp}} .
\ee

The path integral,
\be
Z_{{\perp}} = \int D\Phi_{\perp}
\; \ex{S_{\perp}},
\ee
is essentially the same one that was discussed in the rank 3 case so
it gives a factor of
\be
\left(|{\mathrm H}_{1}(M, {\Bbb Z})|'\right)^{n}.
\ee
It is the first factor that is of interest,
\be
Z_{0} = \int d\Phi_{0} \; \ex{ -S_{0} } .
\ee
The integral $\mu(M) = \int_{M} gd g $, is a well-known invariant of
the manifold $M$ (see \cite{L} for the
Poincar\'{e} dual linking number definition). The Casson invariant in
this case is $-|{\mathrm H}_{1}(M, {\Bbb Z})|' \, \mu(M)$.

We have then that
\be
Z_{X}^{RW}[M] = c'_{X} \left(\lambda_{M}\right)^{n} ,
\ee
where,
\be
c_{X}' = \int d \Phi_{0} \, \frac{\left
( \Lambda_{I}\eps^{IJ}\Lambda_{J}\right)^{n} }{2^{n}n!} .
\ee
Below, we will establish that $c_{X}' = c_{X} $.

\subsection{On the Relationship with [BeHa]}

In \cite{BeHa} the LMO invariant and the Lescop invariant for manifolds with
$b_{1}(M)=2 $ were related. These authors established that the
coefficients of the powers of the Lescop invariant are related to
evaluations of certain diagrams that we can refer to as ${\mathbf H}$
diagrams. One should, as we have previously seen, think of the
vertices in the Rozansky-Witten theory as if they are 3-point
vertices, the $\eta^{I}_{0}$ leg being thought of as the `coupling
constant' (i.e. one focuses on the order of the $\eta^{I}_{0}$ in the
expansion). This 3-point vertex,
\be
V_{IJK}^{3}\, \chi^{I}_{1}\chi^{J}_{2}\chi^{K}_{3} =
\frac{1}{6}\, \Omega_{IJKL}\, \chi^{I}_{\alpha}\chi^{J}_{\beta}
\chi^{K}_{\gamma} \, \eta^{L}_{0} \, \eps^{\alpha \beta \gamma} ,
\label{vert3}
\ee
is what appears finally in (\ref{lefttodo}) in the rank 3
case. Ignoring the $\eta^{I}_{0}$, we see that the vertex carries $3$ legs
which are attached to the three $\chi_{\alpha}^{I}$ zero modes. Each
of the legs carries a different value of $\alpha = 1, 2,3$. In the
current situation, however, we find in the exponent not a 3-point
vertex but rather a 4-point vertex (which is quadratic in the
coupling constant)
\be
V^{4}_{IJKL}\, \chi^{I}_{1}\chi^{J}_{2}\chi^{K}_{1}\chi^{L}_{2} =
\frac{1}{2}\Lambda_{I}\,\eps^{IJ}\Lambda_{J} . \label{vert4}
\ee
This vertex is the ${\mathbf H}$ diagram. The vertex is really a join
of two 3-point vertices along the leg marked 3. The external
legs can only carry the labels $1,2$.

We note that there is another way of expressing the integrals that
still need to be performed and which exhibits very clearly that the
4-point vertex comes from the join of 3-point
vertices. Introduce another Grassmann odd variable
$\psi^{I}$. Then we have\footnote{Setting $a=0$ in (\ref{32}) shows
that the normalization is $\int d\mu(\psi) \, \exp{ ( -1/2
\psi^{I}\eps_{IJ}\psi^{J} )} = 1$, while differentiating twice with
respect to $\Lambda$ and setting $a=1$ and $\Lambda =0$ shows that
$<\psi^{I}\psi^{J}> = \int d\mu(\psi) \, \exp{ ( -1/2
\psi^{K}\eps_{KL}\psi^{L} )} \, \psi^{I}\psi^{J} = \eps^{IJ}$.}
\be
\int d\mu( \psi) \, \ex{\left(- \frac{1}{2}\psi^{I}\eps_{IJ} \psi^{J}
+ a\Lambda_{I}\psi^{I} \right)} = \ex{ \left(-\frac{1}{2} a^{2}
\Lambda_{I} \eps^{IJ} \Lambda_{J} \right)}. \label{32}
\ee
In this way it is as if we have an extra $\chi$ harmonic mode, that is
$\psi^{I}$ plays the role of $\chi_{3}^{I}$. The 4-point
vertex then is really understood as
\be
V^{4}_{IJKL}\, \chi^{I}_{1}\chi^{J}_{2}\chi^{K}_{1}\chi^{L}_{2} =
\frac{1}{2}V^{3}_{I_{1}J_{1}K_{1}}\, \chi^{J_{1}}_{1} \chi^{K_{1}}_{2}\,
\langle  \psi^{I_{1}}\psi^{I_{2}} \rangle \, V^{3}_{I_{2}J_{2}K_{2}}
\, \chi^{J_{2}}_{1} \chi^{K_{2}}_{2} ,
\ee
meaning a contraction of two 3-point vertices along the legs
marked with a 3.

So far we have shown how the ${\mathbf H}$ diagrams arise in the
Rozansky-Witten theory. Now we will see that this characterization of
the ${\mathbf H}$ diagram automatically establishes that the constants
$c_{X}$ for rank 3 manifolds and $c_{X}'$ for rank 2 manifolds are
equal. The expressions (\ref{32}) look like those obtained for the rank 3
case. In fact the resemblance becomes equality with the following
observation: To saturate the integral over $\eta^{I}$ one must expand
the exponential,
\be
\ex{a\Lambda_{I}\psi^{I}} ,
\ee
out to the $2n$'th term. However, in so doing we will also have $2n$
products of $\psi^{I}$, which is exactly what is required to be able to
perform the $\psi^{I}$ integral. Consequently, the only term of the
expansion of the exponential of the quadratic term, $\exp{(-1/2
\psi^{I}\eps_{IJ} \psi^{J} )}$, is the zeroth order piece, namely $1$.

Consequently, we have that
\bea
\int d\Phi_{0}\; \ex{-\frac{1}{2}\mu(M) \Lambda_{I}\eps^{IJ}\Lambda_{J}}& =&
\int d\Phi_{0}\; \ex{-\sqrt{\mu(M)} \Omega_{IJKL} \psi^{I} \chi^{J}_{1}
\chi^{K}_{2} \eta^{L} } \nonumber \\
&=& \mu(M)^{n} \, {\mathbf e}(X) ,
\eea
where the measure on the right hand side of the first equality
includes that of the $\psi^{I}$ field. Now we are done since,
\be
\lambda(M) = {\mathbf e}(X_{AH}) |H_{1}(M, {\Bbb Z})|' \, \mu(M),
\ee
we have shown that
\be
Z_{X}^{RW}[M] = c_{X} \lambda(M)^{n} .
\ee

\section{Manifolds with $b_{1}(M)=1$}\label{manb1}

When $b_{1}(M)=1$, the vertex $V_{1}$ can have at most one harmonic
$\chi^{I}$. Note that in order to saturate the integral over the harmonic
$\chi^{I}$ fields, there is a bound
\be
p+q \geq 2n .
\ee
This bound is already implied by (\ref{p+q}) and so appears to convey no new
information. However, one should read it in a different way. It tells
us that the equality can be met only if one of the $\chi^{I}$ that
appears in the $V_{1}$ vertex and the one that appears in the $V_{2}$
vertex are harmonic. Let $\o$ be a generator for ${\mathrm
H}^{1}(M, {\Bbb Z})$ so that we may write $\chi^{I}_{0} = c^{I} \o$,
where $\o$ satisfies $\int_{M} \o * \o = 1 $.
As before, if a field appears with a $\perp$ subscript then it is
orthogonal to the harmonic modes while, if it has a 
zero subscript then it is understood to be harmonic. We may as well set
\bea
L_{1} & = & \frac{1}{2} \left( g_{ij} \,d\f^{i}_{\perp} * d \f^{j}_{\perp} -
\gamma_{i}^{AK} \gamma_{j}^{BL}\eps_{AB} \Omega_{IJKL} \chi^{I}_{0}
\eta^{J}_{0} \f^{i}_{\perp}*d \f^{j}_{\perp} + 2 \eps_{IJ}
\chi^{I}_{\perp} *d \eta^{J}_{\perp} \right)
\nonumber \\
& = & g_{I \overline{J}} d_{A}\f^{I}_{\perp} * d\overline{\f}^{\overline{J}}_{\perp} +
\eps_{IJ}\chi^{I}_{\perp}*d \eta^{J}_{\perp} \label{l1}\\
L_{2} &=& \frac{1}{2}\left
( \eps_{IJ}(\f_{0})\chi^{I}_{\perp}d\chi^{J}_{\perp} 
+
\Omega_{IJKL} \chi^{I}_{\perp}\chi^{J}_{\perp}
\chi^{K}_{0}\eta^{L}_{0}  \right) ,
\nonumber \\
&=& \frac{1}{2}\eps_{IJ}(\f_{0})\chi^{I}_{\perp}d_{A}\chi^{J}_{\perp}
.  \label{l2}
\eea
The covariant derivative in (\ref{l2}) is defined by
\be
d_{A}\chi^{J} = d\chi^{J} + A^{J}_{\, K} \chi^{K}
\ee
where the ``connection'' is
\bea
A^{I}_{\, J} &=& -\eps^{IM}(\f_{0}) \Omega_{MNJK}(\f_{0}) \chi^{N}_{0}
\eta^{K}_{0} \\ \label{connection}
&=& a^{I}_{\, J} \, \o ,\label{connection1}
\ee
and
\be
a^{I}_{\, J} =
-\eps^{IM}(\f_{0})\Omega_{MNJK}(\f_{0})c^{N}\eta^{K}_{0} .
\ee

The connection is flat; since the tensors that appear in
(\ref{connection}) depend only on the constant maps, and the fields
there are also harmonic, we are assured that $d A^{I}_{\, J} = 0 $.
Furthermore, as $A^{I}_{\, J}$ is proportional to $\o$, we know that
$A^{I}_{\,  J} \wedge A^{J}_{\, K} =0$, so finally
\be
F_{A} =0 .
\ee
Notice that the connection (\ref{connection}) is symmetric when the
labels are both down, $A_{IJ} = \eps_{IK}A^{K}_{\, J} = A_{JI}$, by
virtue of the symmetry properties of $\Omega_{IJKL}$.

\subsection{A Path Integral for Ray-Singer Torsion}\label{pirs}

We now remind the reader of how one formulates the Ray-Singer Torsion
in terms of path integrals. This is a small variant on the formulation
introduced by Schwarz \cite{S}. Let $V$ be a vector bundle over $M$
with a fixed flat connection $A$. One begins with an action
\be
S_{0} = \int_{M} \frac{1}{2}\eps_{IJ} \chi^{I} d_{A} \chi^{J}\label{S0}
\ee
which makes sense for any $M$ with real dimension $4k+1$ with the
$\chi^{I}$ Grassmann even $2k$ forms with values in $V$ or with real
dimension $4k-1$ and
the $\chi^{I}$ are Grassmann odd $2k-1$ forms with values in $V$. As it
stands  this system is not well prescribed since
\be
S_{0}(\chi^{I} + d_{A}\f^{I} ) = S_{0}(\chi^{I}) ,
\ee
that is, the action enjoys a gauge symmetry. In general
$\f^{I}$ will be a form of one degree less than that of
$\chi^{I}$. The symmetry requires that the connection be flat. Hence on the
space of $\chi^{I}$, denoted by ${\mathcal X}$, there is an
action by the gauge group ${\mathcal G}$ given by
\be
\f \left(\chi^{I}\right) = \chi^{I} + d_{A}\f^{I}, \;\;\;\;\; \f \in
{\mathcal G}.
\ee
We do not wish to integrate over ${\mathcal X}$ but rather
over ${\mathcal X}/{\mathcal G}$. Equivalently we chose to
integrate on a slice (section). In doing this one needs to compare
Riemannian volumes on the section and on the space ${\mathcal
X}$. This comparison shows that one must
multiply by a volume factor, namely the Fadeev-Popov
determinant. Ultimately one finds
\be
\int_{{\mathcal X}/{\mathcal G}} \ex{-S_{0}} =
\int_{{\mathcal X}} \ex{-S_{0}} \, \d \left(s \right) \,
\Delta_{{\mathrm FP}} , \label{partitionf}
\ee
where $s$ denotes the section of choice.

Typically one takes the section to be
\be
s^{I} = d_{A}*\chi^{I} , \label{sectiona}
\ee
since, with respect to the metric on $M$, it projects along the direction
of gauge transformations. The path integral with this choice of gauge
is
\be
\int_{{\mathcal X}/{\mathcal G}} \ex{-S_{0}}= \int_{{\mathcal X}}
D\eta^{I} D\f^{I}D\f^{\overline{I}} \ex{-S_{a}}
,\label{pi1}
\ee
where
\be
S_{a} = S_{0} + \int_{M} \eps_{IJ} \eta^{I}d_{A}*\chi^{J} + \int_{M}
\f^{\overline{J}} d_{A}*d_{A} \f^{I} g_{I\overline{J}} . \label{Sa}
\ee
Notice that the path integral over $\eta^{I}$ is there to give back
the delta function constraint onto the section (\ref{sectiona}) while
the integral over $\f^{i}$ reproduces the Fadeev-Popov determinant.

The sum of the actions of (\ref{l1}) and (\ref{l2}),
\be
S_{b} = S_{0} + \int_{M} \eps_{IJ} \eta^{I}d*\chi^{J} + \int_{M}
\f^{\overline{J}} d*d_{A} \f^{I} g_{I\overline{J}} . \label{Sb}
\ee
almost coincides with (\ref{Sa}). A glance at (\ref{Sb}) suggests we are
quantizing the same starting
action (\ref{S0}), but with a different choice of section,
\be
s^{I}= d*\chi^{I} .
\ee
We would expect that the path integral would not depend on which
choice of section we make use of. This is not quite true due to
the presence of zero modes of the fields. The path integral that we
wish to perform, with action (\ref{Sb}), has the condition that we do
not include an integration over the constant $\f^{i}$ and $\eta^{I}$
modes, nor do we integrate over the harmonic part of
$\chi^{I}$. However, in evaluating the path integral which gives
us the Ray-Singer Torsion, no such restriction is made, since one
requires the cohomology of $d_{A}$ to be acyclic (if it is not acyclic
then one explicitly
projects out the harmonic modes of the twisted Laplacian). The path integral
that we want is then not equal to the path integral for the Ray-Singer
Torsion, but rather is equal to the path integral for the Ray-Singer
Torsion divided by the integration over the harmonic modes (of the
usual Laplacian).

The part of the path integral (\ref{pi1}) over the harmonic modes is
\be
\int d\f^{I} d\f^{\overline{J}} d\eta^{I} dc^{I} \, \ex{- \left(
\eta^{I} \eps_{IJ}A^{J}_{\, K} c^{K} -
\f^{\overline{J}}g_{\overline{J} I} A^{I}_{\, K}A^{K}_{\, L} \f^{L}
\right) }= \frac{\det\left(a_{IJ} \right)}{\det
\left(a_{IK} a^{K}_{\, J}\right) } .
\ee
We have found then, for manifolds with $b_{1}(M)=1$, that the path
integral is\footnote{The prefactor of $(1/2\pi)^{2n}$ is part of the
normalization of the path integral.}
\be
Z_{X}^{RW}[M] = \frac{1}{(2\pi)^{2n}}\int_{X} \det(a) \left(\tau_{{\mathrm
RS}}(a)\right)^{1/2} \label{b1rwa}
\ee
where $\tau_{{\mathrm RS}}(a)$ is the Ray-Singer Torsion for the
connection $A$ of a flat ${\mathrm sp}(n)$ bundle over $M$ (in the
$2n$-dimensional representation). Since in the integrand we must pick
a top form, one may move the $2\pi$ factors into a different position,
\be
Z_{X}^{RW}[M] = \int_{X} \det(a/2\pi) \left(\tau_{{\mathrm
RS}}(a/2\pi)\right)^{1/2} .\label{b2rwa}
\ee

\subsection{A BRST Argument}

There is an equivalent way to express the fact that the path integral
that we are interested in yields (\ref{b1rwa}). We start with a BRST
formulation of the model given by the action (\ref{S0}). The BRST
symmetry in question is
\be
Q \chi^{I} = d_{A}\f^{I}, \;\; Q\f^{I} = 0, \;\; Q \f^{\overline{I}} =
g^{\overline{I}J}\eps_{JI} \eta^{I} ,\;\; Q\eta^{\overline{I}} = 0 ,
\label{brst}
\ee
and $Q^{2}=0$. We may still decompose the fields as to whether or not
they are harmonic with respect to the usual de-Rham operator $d$ so
that (\ref{brst}) splits as
\be
\begin{array}{ll}
Q\chi^{I}_{\perp} = d_{A}\f^{I}_{\perp} & Q\chi^{I}_{0} = A^{I}_{\;
J}\f^{J}_{0} \\
Q \f^{\overline{I}}_{\perp} =
g^{\overline{I}J}\eps_{JI} \eta^{I}_{\perp} & Q \f^{\overline{I}}_{0} =
g^{\overline{I}J}\eps_{JI} \eta^{I}_{0} , 
\end{array}
\ee
and
\be
Q\f^{I}_{\perp} = 0 = Q\f^{I}_{0} \;\;\;  Q\eta^{I}_{\perp}= 0=
Q\eta^{I}_{0} .
\ee
Notice that the action (\ref{S0}) is in fact
\be
S_{0} = \int_{M} \frac{1}{2}\eps_{IJ} \chi^{I}_{\perp} d_{A} \label{S02}
\chi^{J}_{\perp} .
\ee
We take the holonomy of $X$ to be irreducible\footnote{This means that
${\mathrm
H}^{0}_{A}(M, V) = 0$. The Ray-Singer torsion (up to some power) is
defined to be the path integral over the fields in (\ref{brst}) with
the proviso that they are orthogonal to the harmonic modes of
$d_{A}$. As we have seen, this means that there are no restrictions on
$\f^{I}$, $\f^{\overline{I}}$ or $\eta^{\overline{I}}$ and the condition on
$\chi^{I}_{\perp}$ is the same for both this theory and the one of
real interest given in (\ref{l1}, \ref{l2}). }, which
means that $A^{I}_{\; K} \f^{K} =0$ implies that $\f^{I} = 0$. This is
not really a restriction since we will not be using any special properties
of the curvature 2-form in any case. So, put another way, we are
considering $sp(n)$ matrices $a^{I}_{\; J}$ which are invertible. The
fact that $\chi_{0}^{I}= \omega c^{I}$ makes no appearance in the
action is due to the fact that it is pure gauge (that is it can be gauge
transformed to zero), since
\be
\omega c^{I} = d_{A} (a^{-1})^{I}_{\; J} c^{J} . \label{pgauge}
\ee

We may
gauge fix in two stages. Firstly we fix the $\perp$ modes, and the gauge
fixing term is taken to be
\be
\{ Q , \int_{M}\, g_{I \overline{J}} \chi^{I}_{\perp} * d
\f^{\overline{J}}_{\perp} \} = \int_{M}\,\left( \eps_{I J} \chi^{I}_{\perp} * d
\eta^{J}_{\perp} + g_{I \overline{J}} d_{A}\f^{I}_{\perp}*d
\f^{\overline{J}}_{\perp}\right)   . \label{gperp}
\ee
Up to this point we see that the path integral that we are interested
in has as its action (\ref{S02}) and (\ref{gperp}). However, the path
integral for the Ray-Singer torsion requires that we also gauge fix
the zero mode symmetry. In fact we should set
$c^{I}$ to zero since, by (\ref{pgauge}) it is pure gauge. In order to
do this we add
\be
\{ Q , \int_{M}\, g_{I \overline{J}} \chi^{I}_{0} * \omega 
\f^{\overline{J}}_{0} \} = \left( \eps_{I
J} c^{I} \eta^{J}_{0} + g_{I \overline{J}} a^{I}_{\; K}\f^{I}_{0} \,
\f^{\overline{J}}_{0}\right) . \label{hamact}
\ee
While this two step gauge fixing is not the usual covariant gauge
fixing we know that, nevertheless, it leads correctly to the
Ray-Singer torsion, since it differs from the covariant gauge fixing
terms by a BRST exact term. Denote the path integral for the
Ray-Singer torsion by $Z_{RS}$, the path integral that one gets simply
by integrating over the perpendicular modes, with action (\ref{S02})
and  (\ref{gperp}), by $Z_{\perp}$, and the zero mode partition
function by $Z_{0}$ that comes from (\ref{hamact}). The path integral
now nicely factors as
\be
Z_{RS}= Z_{\perp}\, Z_{0} ,
\ee
or put another way
\be
Z_{\perp} = Z_{RS}/Z_{0} = \det(a) \left(\tau_{{\mathrm
RS}}(a)\right)^{1/2}.
\ee

\subsection{Explicit Expression for the Path Integral}

Since $a$ takes its values in the adjoint representation of $sp(n)$,
by a global gauge transformation we may rotate it into the Cartan
sub-algebra of $sp(n)$. Denote the $n$ eigenvalues of $a/2\pi$ by $x_{i}$. Then
$a/2\pi$ is conjugate to ${\mathbf diag}\left(x_{1},x_{2},\dots ,x_{n},
-x_{1},  -x_{2}, \dots , -x_{n} \right)$. 

Set $t_{i} = \ex{x_{i}}$. Rewrite the
integrand of (\ref{b2rwa}) as
\be
\det(a/2\pi) \left(\tau_{{\mathrm RS}}(a/2\pi)\right)^{1/2} = \prod_{i=1}^{n}
\, x_{i}^{2} \, \tau_{{\mathrm RF}}(t_{i}) . \label{rstprod}
\ee
In (\ref{rstprod}) $\tau_{{\mathrm RF}}$ is the Reidemeister-Franz
Torsion. The ${\mathrm sp}(1)$ Ray-Singer Torsion for a connection
\be
A = \frac{a}{2\pi} \left( \begin{array}{cc}
1 & 0 \\
0 & -1
\end{array} \right) \o
\ee
is the square of the Reidemeister-Franz Torsion for $a/2\pi$.

\subsection{Reidemeister-Franz Torsion and the Alexander Polynomial}

The relationship between the Reidemeister-Franz Torsion and the
Alexander Polynomial allows us to re-express the Rozansky-Witten
invariant in a form that will prove useful for comparison to known
results about both the Lescop invariant and the LMO invariants. It is
known that the Reidemeister-Franz Torsion and the Alexander
Polynomial\footnote{$\Delta_{M}(t)$ is normalized so as to be
symmetric in $t$ and $t^{-1}$ and so that $\Delta_{M}(1)=1$.}
for a compact closed 3-manifold $M$, $\Delta_{M}(t)$, are related
by \cite{Tu}
\be
\tau_{RF}(M;t) = \frac{\Delta_{M}(t)}{(t^{1/2}-t^{-1/2})^{2}} .
\ee
On substituting $t= \ex{x}$, we see that this relationship may be rewritten as
\be
x^{2}\, \tau_{RF}(M;\ex{x}) =
\left(\frac{x/2}{\sinh{x/2}}\right)^{2}\Delta_{M}(\ex{x}) ,
\ee
so that
\be
\prod_{i=1}^{n}x^{2}_{i}\, \tau_{RF}(M;\ex{x_{i}}) = \hat{{\mathrm A}}(X)
\prod_{i=1}^{n} \Delta_{M}(\ex{x_{i}}) .
\ee
We have then\footnote{The sign has been fixed in \cite{RW} and \cite{T}.}
\be
Z_{X}^{RW}[M] = -\int_{X} \; \hat{{\mathrm A}}(X)
\prod_{i=1}^{n} \Delta_{M}(\ex{x_{i}}) .
\ee

Given any compact hyper-K\"{a}hler manifold $X$ of dimension $4n$ the
Todd genus is $n+1$, so that (\cite{RW}, \cite{T})
\be
Z_{X}^{RW}[S^{2}\times S^{1}] = -(n+1) ,
\ee
since $\Delta_{S^{2}\times S^{1}}(t)=1$. The
$S^{[n]}$ series is then enough to guarantee that there is a non-zero
invariant in every dimension for $S^{2}\times S^{1}$.

Since the Chern characters of the tangent bundle of $S^{[n]}$
have been shown to be non-zero up to $n=8$, we know that all the
invariants are realized up to this degree for any 3-manifold $M$ with
$b_{1}(M)=1$.

\subsection{On the Relationship with [GH]}

In this section we assume some familiarity with the notation used in
\cite{GH}. 
According to \cite{GH} for rank 1 manifolds with no torsion in ${\mathrm
H}_{1}(M, {\Bbb Z})$ ($={\Bbb Z}$), the LMO invariant may be written
as
\be
Z^{LMO}(M) = \langle \exp{ _{\sqcup} \alpha (M)} \rangle . \label{LMO1}
\ee
The notation is as follows; $\alpha (M)$ corresponds to a particular
set of diagrams (more on this below), the cup $_{\sqcup}$ means take the
disjoint union of diagrams (the exponential is to be understood in the
same way) the brackets $\langle \; \rangle$ mean contraction over all
external legs in all possible ways.

We now need to explain what $\alpha(M)$ is. Let
\be
-\frac{1}{2} \log{ \left( \Delta_{M}( \ex{x}) \right) } =
\sum_{m=1}^{\infty} a'_{2m}(M) \, x^{2m} , \label{lalex}
\ee
and
\be
- \log{ \left( \frac{x/2}{\sinh{(x/2)}} \right)} = \sum_{m=1}^{\infty}
2b_{2m} \, x^{2m} .
\ee
The logarithm of our favourite product is then
\be
-\frac{1}{2}\log{\left( \left( \frac{x/2}{\sinh{(x/2)}} \right)^{2}
\Delta_{M}(\ex{x}) \right)} = \sum_{m=1}^{\infty} \left(2b_{2m}+
a'_{2m}(M) \right) \, x^{2m} .
\ee
By definition we have
\be
\alpha(M) = \sum_{m=1}^{\infty} \left(2b_{2m}+
a'_{2m}(M) \right) \, \omega_{2m} ,
\ee
where the powers of the indeterminant $x$ have been replaced by the
wheel diagrams $\omega_{2m}$.

How does this data compare with what we have just learnt about the
Rozansky-Witten invariant? The integrand is expressible as
\be
\exp{\left(-2 \sum_{m=1}^{\infty} \left(2b_{2m}+
a'_{2m}(M) \right) \, .\Tr \left(\frac{R}{2\pi}\right)^{2m}\right) }
\ee
The correspondence is thus quite clear. The diagrams $\omega_{2m}$ are
replaced by $-2\Tr \left(\frac{R}{2\pi}\right)^{2m}$, the cup product
should be read as the wedge product and contraction over all legs becomes
integration over $X$.  Under the conditions that they consider,
$b_{1}(M)=1$ and $|
{\mathrm H}_{1}(M, {\Bbb Z})|' = 1 $
the LMO invariant is determined by and indeed determines the Alexander
Polynomial. The same is true for the Rozansky-Witten invariants (we do
not require that $| {\mathrm H}_{1}(M, {\Bbb Z})|' = 1$). Clearly the
Alexander Polynomial determines $Z^{RW}_{X}[M]$. The converse is also
true, since, as one
increases the dimension of $X$, higher and higher derivatives of the
Alexander Polynomial make an appearance. When, $\dim_{{\Bbb R}}X=4n$,
one has $\Delta^{(2n)}_{M}(1)$ and lower order derivatives appearing
on the right hand side and so one may, inductively, determine the
Taylor series of the Alexander polynomial around $t=1$.

\section{The Hilbert Space on $S^{2}$}\label{manb10}

The techniques that we have been using thus far do not suffice to give
closed form expressions for the invariants when the rank is zero. We
need to make use of another point of view on the path
integral.\footnote{This may well be related to the fact that the
invariants are truly effective only for $b_{1}(M)=0$.}

A path integral on a manifold with boundary prepares (meaning is) a
vector in a Hilbert space of states. More often than not the Hilbert
space is really an infinite dimensional Fock space. In some very
special circumstances the `Hilbert' space is a finite dimensional
vector space and has some properties that we would like it to have,
for example it comes equipped with a non-degenerate inner product. The
inner product may not be (and here is not) positive, but the vector space
will be called a Hilbert space of states. It was conjectured in
\cite{RW} that providing one chooses the
hyper-K\"{a}hler manifold $X$ to be compact, fortune smiles on us and
the Hilbert spaces of states for the Rozansky-Witten theory are related
to certain cohomology groups of $X$ (so that, in particular, they are finite
dimensional). Supporting evidence for this was provided in \cite{T}
and we will use this fact below.

\subsection{The Hilbert Space for (non) Compact $X$}

If a manifold $X$ is non-compact one has a
choice of which cohomology one considers for the space. For example 
when $X_{{\mathrm SU}(2)} = X_{AH}$, by the work of Sen \cite{Se} on
the duality
conjecture of $N=4$ super Yang-Mills theory in four dimensions, it is
known that in ${\mathrm L}^{2}$-cohomology there is only one non-zero
form, which must be a $(1,1)$ form. On the other
hand, when we come to consider the properties of the invariants under
connected sum we
will see that the total number of states available depends on the
cohomology groups ${\mathrm H}^{(0,p)}$, which we would say is empty
for $X_{AH}$ if they are understood as ${\mathrm H}^{(0,p)}_{{\mathrm
L}^{2}} $. This would have the nasty consequence that the
Rozansky-Witten invariant would vanish for $S^{2} \times S^{1}$. 

Different choices of cohomology for non-compact $X$ can,
therefore, lead to wildly different results. It is for this reason
that we make use, in the
following, of compact hyper-K\"{a}hler manifolds. It was argued in
\cite{RW} that there is no loss in doing so since, in any case, the
only dependence on $X$ is through integrals like (\ref{inttypes}).
Thus if we take the manifolds to be compact, we can consider the usual
cohomology theory of these spaces. Once the dependence on $X$ is
worked out for arbitrary compact $X$ and written in terms of integrals of the
form (\ref{inttypes}) then this dependence will be correct also for
$X_{G}$. 

However, in the case of non-compact manifolds one needs to make use of
a slightly different set of basic observables \cite{RW}. For example,
they may be based on powers of
\be
{\cal O} = \left( \eps^{J_{1}J_{2}}\eps^{K_{1}K_{2}}
\eps^{L_{1}L_{2}} \, \Omega_{I_{1}K_{1}J_{1}L_{1}} \Omega_{I_{2} K_{2}
J_{2}L_{2}}\right) \eta^{I_{1}}\eta^{I_{2}} ,
\ee
suitably normalized, or on other combinations of the curvature and
holomorphic two form. There are a number of conditions that must be
met by these operators. Firstly, whatever these operators are, in the
compact case they must be equivalent (cohomologous) to
the original observables ${\mathcal O}(k)$. Furthermore insertions of
these operators in the path integral ought to lead to integrable expressions
(on $X$). Finally, they should obey the condition (\ref{norm}).

\subsection{Partition Function in terms of the Hilbert Space}

The boundary in the present setting is a Riemann surface
$\Sigma$. Given such vectors in the Hilbert space one can use the usual tenets
of quantum field theory to reconstruct the partition function
$Z_{X}^{RW}[M]$. For example the path integral of the field theory on
the 3-ball $B^{3}$ prepares a
state on the boundary $S^{2}$ (the Hilbert space will be described
shortly). We denote this state (vector) by
\be
| \psi^{(0)}> = | B^{3} > . \label{b1}
\ee
Let $B^{3} \subset M$ a 3-ball
inside $M$. Then the path integral on $M\bs B^{3}$ will give us
another vector in the Hilbert space which is denoted by
\be
\mid M\bs B^{3} > .
\ee
The path integral tells us that the partition function on $M$ is given as the
inner product of the two vectors
\be
Z_{X}^{RW}[M] = <B^{3} \mid M\bs B^{3} >.
\ee
In this notation one has that if $M$ is a connected sum $M=M_{1}\#
M_{2}$, where
\be
\partial M_{1} = S^{2} = \partial M_{2} ,
\ee
then
\be
Z_{X}^{RW}[M] = <M_{2}^{*} \mid M_{1} > =
<M_{1}^{*}\mid M_{2} > ,
\ee
where a star superscript means take the opposite orientation. There is
a more general formula available, which is obtained by similar
arguments. Let $M$ have a Heegaard decomposition along a genus $g$
Riemann surface $\Sigma_{g}$ as $M= M_{1}\#_{\Sigma_{g}}M_{2}$, where
\be
\partial M_{1} = \Sigma_{g} = \partial M_{2} ,
\ee
then
\be
Z_{X}^{RW}[M] = <M_{2}^{*} \mid M_{1}> =
<M_{1}^{*} \mid M_{2} > .
\ee

We will make use of formulae of the above type to establish the
properties of the generalized Casson invariant under the operation of
connected sum. Before that we need a digression on the invariants of
$S^{3}$.

\subsection{The Path Integral on $S^{3}$}
Rozansky and Witten have established that their generalized invariant
$Z_{X}^{RW}[M]$ under change of orientation behaves as
\be
Z_{X}^{RW}[M^{*}] = (-1)^{n(1+b_{1}(M))} Z_{X}^{RW}[M] .
\ee
Insertion of the operator ${\mathcal O}(k)$ essentially lowers the
effective dimension of $X$ by $4k$ to $4(n-k)$, so that under
orientation reversal
\be
Z_{X}^{RW}[M^{*}, {\mathcal O}(k)] = (-1)^{(n-k)(1+b_{1}(M))}
Z_{X}^{RW}[M, {\mathcal O}(k)] . \label{orev}
\ee
One may obtain (\ref{orev}) as follows. The orientation properties
are determined by counting the number of $\eps_{\mu \nu \rho}$ tensors that
appear in the perturbative diagrams. There is one such tensor for each
$V_{1}$ vertex and one for each $\chi$ propagator, $<\chi \chi>$, so
the behaviour
under sign reversal is multiplication by
\be
(-1)^{\# V_{1} + \# <\chi \chi>} . \label{osign}
\ee
For a diagram with $p$ $V_{1}$ vertices, $q$ $V_{2}$ vertices and an
insertion of ${\mathcal O}(k)$ we have
\be
p + q + 2k = 2n
\ee
by counting $\eta$ harmonic modes. On the other hand, as discussed
above, when there are $\chi$ zero modes they must appear in the
vertices, so the number of $\chi$'s which are not zero modes in
$V_{1}$ is $3-b_{1}(M)$ and the number in $V_{2}$ is $1-b_{1}(M)$. The
number of $\chi$ propagators is simply one half of the total number of
non-harmonic $\chi$ legs in the diagram and since the legs can only
come from vertices we have
\be
\# <\chi \chi> = \frac{1}{2}\left( (3-b_{1}(M))p + (1-b_{1}(M))q \right) .
\ee
The sign is then determined by
\be
p+ \# <\chi \chi> = \frac{1}{2}(1-b_{1}(M))(p+q) = (1-b_{1}(M))(n-k) .
\ee

If there are orientation reversing diffeomorphisms available on $M$, then
(\ref{orev}) tells us that
\be
Z_{X}^{RW}[M, {\mathcal O}(k)] = (-1)^{(n-k)(1+b_{1}(M))}
Z_{X}^{RW}[M, {\mathcal O}(k)] ,
\ee
which provides a vanishing theorem for some of the invariants on
$M$. For example $S^{3}$ admits an orientation reversing
diffeomorphism, so that the Casson invariant (corresponding here to
$n=1$, $k=0$) enjoys
\be
\lambda (S^{3}) = - \lambda (S^{3}) = 0 .
\ee
In fact from (\ref{orev}) we learn that for manifolds with $b_{1}(M)=0$
which admit orientation reversing diffeomorphisms, the invariants
$Z_{X}^{RW}[M, {\mathcal O}(k)] $ will necessarily vanish unless $n-k
= 2m$ for some $m$.

We would have liked the slightly stronger result that
$Z_{X}^{RW}[S^{3},{\mathcal O}(k)] $ will necessarily vanish unless
$n-k = 0$, this, however, does not seem to be available. Potentially
this is worrisome as one would expect that a generalized Casson invariant for
$M$ such that $\pi_{1}(M)=1$ would vanish. (For example the $SU(3)$
invariant of Boden and Herald \cite{BH} is designed
to vanish for $M$ such that $\pi_{1}(M)=1$.) We suggest, therefore,
that the generalization of the Casson invariant for a rational
homology sphere that matches the gauge theoretic $SU(3)$ one is not
$Z_{X_{SU(3)}}^{RW}[M]$, but rather
\be
\lambda_{SU(3)}(M) = \lambda_{X_{SU(3)}}(M) = Z_{X_{SU(3)}}^{RW}[M] -
|{\mathrm H}_{1}(M, {\Bbb Z})|
\, Z_{X_{SU(3)}}^{RW}[S^{3}] , \label{cass3}
\ee
where $X_{SU(3)}$ is the reduced $SU(2)$ $3$ monopole moduli space. Not only
does this satisfy the requirement that $\lambda_{X_{SU(3)}}(S^{3})=0$,
but it also has good properties under connected sum, as we will
see. For $\ZHS$'s our proposal amounts to
\be
\lambda_{SU(3)}(M) = Z_{X}^{RW}[M] -  Z_{X}^{RW}[S^{3}] .
\ee

\subsection{The Hilbert Space on $S^{2}$ and the Connected Sum Formula}

The Hilbert space of states on a Riemann surface was described in
\cite{RW}. Here we will look at the small Hilbert space of
states, ${\mathcal H}_{\Sigma}$, those states which are $Q$ invariant
(modulo exact terms). The $Q$ operator is identified as the Dolbeault
operator $\overline{\partial}$ on $X$. The small Hilbert space will be
related to $\overline{\partial}$-cohomology of certain classes of
forms on $X$.

The Hilbert space of states on $S^{2}$, ${\mathcal H}_{S^{2}}$, for a
compact hyper-K\"{a}hler manifold $X$ is
\be
{\mathcal H}_{S^{2}} = \bigoplus_{k=0}^{2n} {\mathrm H}^{(0,k)}(X) .
\ee
The following result, which follows from Berger's classification
theorem on the holonomy of Riemannian manifolds, is useful (for this and some
other useful information on
hyper-K\"{a}hler manifolds one may consult \cite{Be}). The holonomy
group of $X$ is a subgroup of
${\mathrm sp}(n)$ and if $X$ is irreducible, then it is ${\mathrm
sp}(n)$. Let ${\mathrm h}^{(p,q)} = \dim {\mathrm H}^{(p,q)}$, for
irreducible $X$,
\bea
{\mathrm h}^{(0,k)}(X) &=& 0, \; \forall \, k \; {\mathrm odd} \nonumber \\
{\mathrm h}^{(0,k)}(X) &=& 1, \; \forall \, k \; {\mathrm even}.
\eea
This means that the real dimension of ${\mathcal H}_{S^{2}}$ is
$n+1$. Furthermore, the elements in ${\mathrm H}^{(0,2k)}(X)$ are
generated by the $k$-th exterior power of the holomorphic symplectic
2-form $\eps$. Consequently any vector in ${\mathcal H}_{S^{2}}$ can
be expressed as
\be
v = \oplus_{k=0}^{n} v_{k} \eps^{k} .
\ee

Denote the path integral that includes insertion of the observable
corresponding to $\eps^{k}$ in the 3-ball by
\be
| \psi^{(2k)}> = | B^{3}, {\mathcal O}(k) > . \label{bk}
\ee
Let the states
\be
<\psi^{(2k)} | = < ( B^{3})^{*}, {\mathcal O}(k)| \, ,
\ee
be defined in such a way that,
\be
Z_{X}^{RW}[S^{3}, {\mathcal O}(p+q)] = < ( B^{3})^{*},
{\mathcal O}(p)| B^{3}, {\mathcal O}(q) > .
\ee

Let $M$ be a rational homology sphere and $B^{3} \subset M$ a 3-ball
inside $M$. Since ${\mathcal H}_{S^{2}}$ is $n+1$-dimensional, the path
integral on $M \bs B^{3}$ yields a state that can be expanded in the basis
generated by (\ref{bk}). This state is then
\be
| M\bs B^{3} > = \sum_{k=0}^{n} \, \lambda^{k}_{X}(M) \, |
\psi^{(2k)}> , \label{mb3}
\ee
for some coefficients $\lambda^{k}_{X}(M)$.

To determine the coefficients in (\ref{mb3}) we
require some properties of the path integral on $S^{3}$. Before
proceeding to the general case we review the way that Rozansky and Witten
derived the connected sum formula for the Casson invariant, $n=1$, and
then derive the analogous expressions for $n=2$.

\subsubsection{$\dim_{{\Bbb R}}X=4$}
When $n=1$, the
Hilbert space is 2-dimensional and (\ref{mb3}) is simply
\be
| M\bs B^{3} > = \lambda_{X}^{0}(M) \, | \psi^{(0)}> + \lambda^{1}_{X}(M)
\, | \psi^{(2)}>.
\label{mb4}
\ee
In this case we have the inner products
\be
Z_{X}^{RW}[S^{3}, {\mathcal O}(k)] = \left\{ \begin{array}{cl}
 1 & {\mathrm for}\;\; k=1 \\
 0 & {\mathrm for}\;\; k=0
\end{array}
\right. . \label{s3}
\ee
We deduce that
\bea
Z_{X}^{RW}[M, {\mathcal O}(k)] &=&  \lambda_{X}^{0}(M) \, <\psi^{(2k)} |
\psi^{(0)}> + \lambda^{1}_{X}(M)
\, <\psi^{(2k)} | \psi^{(2)}> \nonumber \\
& = & \lambda_{X}^{0}(M)Z_{X}^{RW}[S^{3}, {\mathcal O}(k)] +
\lambda^{1}_{X}(M) Z_{X}^{RW}[S^{3}, {\mathcal O}(k+1)] \nonumber \\
& = & \lambda_{X}^{0}(M) \d_{k,1} + \lambda^{1}_{X}(M)\d_{k,0} ,
\eea
or put another way
\bea
\lambda^{0}_{X}(M)& = & Z_{X}^{RW}[M, {\mathcal O}(1)] = |H_{1}(M,
{\Bbb Z})| ,\nonumber \\
\lambda_{X}^{1}(M) & = & Z_{X}^{RW}[M] ,
\eea
and
\be
| M\bs B^{3} > =  Z_{X}^{RW}[M]\, | \psi^{(2)}> + |H_{1}(M,{\Bbb Z})|
\, | \psi^{(0)}>.
\ee
One can likewise ascertain that
\be
| M\bs B^{3}, {\mathcal O}(1) > = |H_{1}(M,{\Bbb Z})| \, | \psi^{(2)}> .
\ee
The state $<(M\bs B^{3})^{*}|$ designates the path integral on the
manifold with boundary $S^{2}$ but with opposite orientation. One can
expand this state also as
\be
<(M\bs B^{3})^{*}| = Z_{X}^{RW}[M]\, <\psi^{(2)}| + |H_{1}(M,{\Bbb Z})|
\, < \psi^{(0)}|
\ee
so that
\bea
& & Z_{X}^{RW}[M_{1}\#M_{2}]  =   |H_{1}(M_{1},{\Bbb Z}) |Z_{X}^{RW}[M_{2}]+
|H_{1}(M_{2},{\Bbb Z})|  Z_{X}^{RW}[M_{1}] \label{con2a} \\
& & Z_{X}^{RW}[M_{1}\#M_{2},{\mathcal O}(1) ]  =  |H_{1}(M_{1},{\Bbb Z})
| |H_{1}(M_{2},{\Bbb Z})| = |H_{1}(M_{1}\#M_{2},{\Bbb Z})| . \label{con2b}
\eea
These formulae tell us all we need to know about the connected sum
properties of the Casson invariant.

\subsubsection{$\dim_{{\Bbb R}}X=8$}
In this section we will repeat the calculation of the behaviour of the
path integral under connected sum for $n=2$. This time the Hilbert
space is 3-dimensional, so for (\ref{mb3}) we have
\be
| M\bs B^{3} > = \lambda_{X}^{0}(M) \, | \psi^{(0)}> + \lambda^{1}_{X}(M)
\, | \psi^{(2)}> + \lambda^{2}_{X}(M) | \psi^{(4)}> .
\ee
The problem that we face has to do with the inner product. This time
we only know that
\be
Z_{X}^{RW}[S^{3}, {\mathcal O}(k)] = \left\{ \begin{array}{cl}
 1 & {\mathrm for}\;\; k=2 \\
 0 & {\mathrm for}\;\; k=1 \\
 a & {\mathrm for}\;\; k=0
\end{array}
\right. . \label{s38}
\ee
We have not determined the constant $a =Z_{X}^{RW}[S^{3}] $ and see no
obvious way of specifying it,\footnote{In this case one cannot appeal to
orientation reversal to rule it out.} without, that is, resorting to
a direct calculation or by appealing to Chern-Simons theory. However,
let us see how far we can go without knowing $a$.

We have
\bea
& & Z_{X}^{RW}[M, {\mathcal O}(k)] \nonumber \\
& & \; = \lambda_{X}^{0}(M) \, <\psi^{(2k)} |
\psi^{(0)}> + \lambda^{1}_{X}(M)
\, <\psi^{(2k)} | \psi^{(2)}> + \lambda^{2}_{X}(M) <\psi^{(2k)}|
\psi^{(4)}> \nonumber \\
& & \;  =  \lambda_{X}^{0}(M)Z_{X}^{RW}[S^{3}, {\mathcal O}(k)] +
\lambda^{1}_{X}(M) Z_{X}^{RW}[S^{3}, {\mathcal O}(k+1)]  \nonumber \\
& & \;\;\;\; +
\lambda^{2}_{X}(M) Z_{X}^{RW}[S^{3}, {\mathcal O}(k+2)] \nonumber \\
& &\; =  \lambda_{X}^{0}(M) (a\d_{k,0} + \d_{k,2}) +
\lambda^{1}_{X}(M) \d_{k,1} +  \lambda^{2}_{X}(M) \d_{k,0} .
\eea
We deduce therefore that
\bea
\lambda_{X}^{0}(M) &=&  Z_{X}^{RW}[M, {\mathcal O}(2)] = |H_{1}(M,{\Bbb
Z})|^{2} \nonumber \\
\lambda^{1}_{X}(M) & = & Z_{X}^{RW}[M, {\mathcal O}(1)] \nonumber \\
\lambda^{2}_{X}(M) & = & Z_{X}^{RW}[M] - a |H_{1}(M,{\Bbb Z})|^{2} .
\eea
As previously explained, the insertion of the operator ${\mathcal
O}(1)$ reduces the effective dimension of the hyper-K\"{a}hler
manifold $X$ by $4$, so that $Z_{X}^{RW}[M, {\mathcal O}(1)]$ is
proportional to the path integral $Z_{K3}^{RW}[M]$. This means that
$\lambda^{1}_{X}(M)$ is proportional to the $SU(2)$ Casson
invariant. We would like to call $\lambda^{2}_{X}(M)$ the degree $2$
Casson invariant.

How do these behave under connected sum? First we note that
\be
<(M\bs B^{3})^{*}| = \lambda_{X}^{0}(M) \, < \psi^{(0)}| + \lambda^{1}_{X}(M)
\, <\psi^{(2)}| + \lambda^{2}_{X}(M) < \psi^{(4)}| ,
\ee
so that
\be
Z_{X}^{RW}[M_{1}\#M_{2}] = <(M_{1}\bs B^{3})^{*}|M_{2}\bs B^{3}>.
\ee
We find that
\bea
Z_{X}^{RW}[M_{1}\#M_{2}]
& = & \lambda_{X}^{0}(M_{1})
\lambda^{2}_{X}(M_{2}) +  \lambda_{X}^{0}(M_{2})\lambda^{2}_{X}(M_{1})
\nonumber \\
& & + \lambda^{1}_{X}(M_{1})\lambda^{1}_{X}(M_{2}) + a
\lambda_{X}^{0}(M_{1})\lambda_{X}^{0}(M_{2}) ,
\eea
which can be put in the nicer form
\be
\lambda^{2}_{X}(M_{1}\#M_{2}) = \lambda_{X}^{0}(M_{1})
\lambda^{2}_{X}(M_{2}) + \lambda^{1}_{X}(M_{1})\lambda^{1}_{X}(M_{2})+
\lambda^{2}_{X}(M_{1})\lambda_{X}^{0}(M_{2}) .
\ee

\subsubsection{$\dim_{{\Bbb R}}X=4n$}
In this section we will provide a proof that the invariants
$\lambda_{X}^{p}(M)$ satisfy a pleasing property under connected sum,
namely (\ref{consumf}) below. In order to do so we will make use of a
property of the inner product in the
basis (\ref{bk}). Let,
\be
<\psi^{(2k)}|\psi^{(2l)}>= G_{k,l}. \label{inprod}
\ee
~From the definitions we have that $G$ is a $(n+1)\times (n+1)$
symmetric matrix where $k,l \;= 0, 1, \dots , n$. G has all unit
entries on the anti diagonal, that is
$G_{k,l}=1$ if $ k+l =n$, and only zero entries below the anti
diagonal, $G_{k,l} = 0$ if $k+l > n$. Consequently $\det G = \pm 1$,
so that in particular G is invertible. All of these properties follow
from the fact that
\bea
G_{k,l}= <\psi^{(2k)}|\psi^{(2l)}> & = & <B^{3}, {\mathcal O}(k)  | B^{3},
{\mathcal O}(l) > \nonumber \\
& = & Z_{X}^{RW}[S^{3}, {\mathcal O}(k+l)] .
\eea

With this inner product in place we now claim that
\be
| M\bs B^{3}, {\mathcal O}(p) > = \sum_{k=0}^{n-p}
\lambda_{X}^{k}(M)|\psi^{(2k+2p)}> . \label{mbo}
\ee
The advantage of having such a formula is that it involves a smaller
number of vectors on the right hand side and so effectively decreases
the size of the Hilbert space that we have to work with. To establish
(\ref{mbo}) we note that quite generally,
\be
| M\bs B^{3}, {\mathcal O}(p) > = \sum_{k=0}^{n-p}
\gamma_{X}^{k}(M)|\psi^{(2k+2p)}> , \label{mbog}
\ee
for some, to be determined, coefficients $\gamma_{X}^{k}(M)$. By construction
\be
<(B^{3})^{*}, {\mathcal O}(q) | M\bs B^{3}, {\mathcal O}(p) > =
<(B^{3})^{*} ,{\mathcal
O}(p+q)  | M\bs B^{3}>,
\ee
which means that
\be
\sum_{k=0}^{n}\gamma_{X}^{k}(M) G_{q,k+p} =
\sum_{k=0}^{n}\lambda_{X}^{k}(M) G_{p+q,k}.
\ee
Taken together with the fact that $G_{q,k+p}= G_{p+q,k}$ and that $G$ is
invertible one has $\gamma_{X}^{k}(M) = \lambda_{X}^{k}(M)$ as
claimed. Thus we have that,
\be
\lambda_{X}^{n-p}(M)= Z_{X}^{RW}[M, {\mathcal O}(p)] -
\sum_{k=0}^{n-p-1}\lambda_{X}^{k}(M)G_{0,k+p} .
\label{ldef}
\ee
This equation is recursive, meaning that the right hand side only
involves $\lambda_{X}^{k}(M)$'s for $k$'s of lower order.

The connected sum formula that we wish to prove is
\be
\lambda_{X}^{p}(M_{1}\# M_{2}) = \sum_{ k+l
=p}\lambda_{X}^{k}(M_{1})\lambda_{X}^{l}(M_{2}) . \label{consumf}
\ee
We will establish (\ref{consumf}) by induction on $p$. Since
$\lambda_{X}^{0}(M) = |{\mathrm H}_{1}(M, {\Bbb Z})|$ (\ref{consumf}) holds for $p=0$.

The connected sum, property is
\bea
& & Z_{X}^{RW}[M_{1}\#M_{2}, {\mathcal O}(n-p-1)] \nonumber \\
& & \;\;\; =  <(M_{1}\bs
B^{3})^{*}|M_{2}\bs B^{3}, {\mathcal O}(n-p-1)>
\nonumber \\
& & \;\;\; = \sum_{k,l =0}^{n} \lambda_{X}^{k}(M_{1})\lambda_{X}^{l}(M_{2})
G_{k,l+n-p-1} \nonumber \\
& & \;\;\; =  \sum_{k+l
=p+1}\lambda_{X}^{k}(M_{1})\lambda_{X}^{l}(M_{2}) \nonumber \\
& & \;\;\;\; \;\;\;\; +
\sum_{0 \leq k+l \leq p} \lambda_{X}^{k}(M_{1})\lambda_{X}^{l}(M_{2})
G_{k,l+n-p-1} . \label{consum}
\eea
If we can show that
\bea
\lambda_{X}^{p+1}(M_{1}\# M_{2}) &=& Z_{X}^{RW}[M_{1}\#M_{2},{\mathcal
O}(n-p-1)  ] \nonumber \\
& & \;\;\; - \sum_{
0\leq k+l  \leq p} \lambda_{X}^{k}(M_{1})\lambda_{X}^{l}(M_{2})
G_{k,l+n-p-1} \label{2show}
\ee
then we will have established (\ref{consumf}). By (\ref{ldef})
\be
\lambda_{X}^{p+1}(M)= Z_{X}^{RW}[M, {\mathcal O}(n-p-1)] -
\sum_{m=0}^{p}\lambda_{X}^{m}(M)G_{0,m+n-p-1} .
\label{ldef2}
\ee
By the inductive hypothesis, $\lambda_{X}^{m}(M)$ satisfies
(\ref{consumf}) for all $k \leq p$, hence
\bea
& &  \sum_{m=0}^{p}\lambda_{X}^{m}(M_{1}\#
M_{2})G_{0,m+n-p-1} \nonumber \\
&& \;\; =   \sum_{m=0}^{p}\sum_{ k+l
=m}\lambda_{X}^{k}(M_{1})\lambda_{X}^{l}(M_{2})G_{0,m+n-p-1}
\nonumber \\
& & \;\; =  \sum_{m=0}^{p}\sum_{ k+l
=m}\lambda_{X}^{k}(M_{1})\lambda_{X}^{l}(M_{2})G_{0,k+l+n-p-1} .
\eea
As the summand in the last line does not depend on $m$, we have $0
\leq k+l \leq p$, whence
\be
\sum_{m=0}^{p}\lambda_{X}^{m}(M_{1}\#
M_{2})G_{0,m+n-p-1} = \sum_{k+l \leq
p}\lambda_{X}^{k}(M_{1})\lambda_{X}^{l}(M_{2})G_{0,k+l+n-p-1}
. \label{css}
\ee
We are done, since for $M= M_{1}\# M_{2}$ (\ref{css}) and (\ref{ldef2}) imply
(\ref{2show}).

Let us state the result of this section again
\be
\lambda_{X}^{p}(M_{1}\# M_{2}) = \sum_{ k+l
=p}\lambda_{X}^{k}(M_{1})\lambda_{X}^{l}(M_{2}) . \label{CONSUM}
\ee
Of course if one would prefer the connected sum formulae for
$Z_{X}^{RW}[M, {\mathcal O}(p)]$ then it is a simple matter to pass to
those given the ones for the $\lambda_{X}^{p}(M)$.

\section{The RW and LMO Invariants}

In order to make contact with the LMO and generalized Casson
invariants we will
find that we need to use normalized Rozansky-Witten invariants. In
both cases, as we will see, the required normalization is such that
the invariant vanishes on $S^{3}$. The reason this is needed to make
contact with \cite{LMO} is basically a question of normalization of
the LMO invariant and so appears here as a direct question of
normalization. On the other hand for the Casson invariant there is a
subtlety which arises in the gauge theoretic setting that requires us
to normalize the invariants in precisely the same way as for the LMO
invariants. The arguments we present for the exact form of the
normalization in the Casson case are suggestive but not complete even
at the physical level of rigour. 

\subsection{Weights in the RW Theory}

To simplify notation, we will write $X_{n}$ for an arbitrary
hyper-K\"{a}hler $4n$ dimensional manifold.

The Rozansky-Witten partition function can be seen to be given by
(c.f. \cite{RW} (3.41)- (3.43))
\be
Z_{X_{n}}^{RW}[M] = \left( | {\mathrm H}_{1}(M, {\Bbb Z} | ' \right)^{n}
\sum_{\Gamma_{n}} b_{\Gamma_{n}}(X_{n}) I_{\Gamma_{n}}(M) ,
\ee
where $\Gamma_{n}$ are Feynman graphs associated with the fact that
the manifold $X$ has dimension $4n$ and the weights 
\be
b_{\Gamma_{n}}(X_{n})= \frac{1}{(2\pi)^{n}} \int_{X} W_{\Gamma}(X,
\f_{0}) \sqrt{g} d^{n}\f_{0}.
\ee
$W_{\Gamma}$ is a product of the tensors $\Omega$ with their indices
contracted by the tensor $\eps^{IJ}$ contained in the $\chi^{I}$
propagator and by $\eps^{I_{1}, \dots, I_{2n}}$ contained in the
$\eta^{I}$ zero-mode expectation value. 

Rozansky and Witten have established
that the weights satisfy the IHX relations. More
generally we will show in this section that
\be
Z_{X_{n}}^{RW}[M, {\mathcal O}(n-k)] = \left( | {\mathrm H}_{1}(M, {\Bbb Z} | ' \right)^{k}
\sum_{\Gamma_{k}} b_{\Gamma_{k}}(X_{n}) I_{\Gamma_{k}}(M) ,
\ee
where $b_{\Gamma_{k}}(X_{n})$ is determined below. This shows that the
insertion of the operators ${\mathcal O}(n-k)$ effectively lowers the
dimension of $X$. The proof of IHX
for $b_{\Gamma_{n}}(X_{n})$, which is essentially a consequence of the
Bianchi identity for $\Omega_{IJKL}$, extends to
$b_{\Gamma_{k}}(X_{n})$. 

Introduce,
\be
Z_{n}^{RW}[M]  = \left( | {\mathrm H}_{1}(M, {\Bbb Z} | ' \right)^{n}
\sum_{\Gamma_{n}} b_{\Gamma_{n}} I_{\Gamma_{n}}(M)
\ee
which, given a hyper-K\"{a}hler manifold $X$, is a map from diagrams
to numbers such that when $X$ is $4n$ dimensional we have
\be
Z_{X}^{RW}[M] = Z_{n}^{RW}[M](X) .
\ee
We need to be more precise. So far we have defined the weights
$b_{\Gamma_{n}}$ as they act on manifolds $X$ of dimension $4n$, we
will need to give a more general definition for their action on
hyper-K\"{a}hler manifolds of arbitrary dimension.

We will now show that insertions of the operators ${\mathcal O}(k)$ 
(where $M$ is a $\QHS$ since otherwise the path integral vanishes), 
have the effect of lowering the effective dimension of $X$. This means
that the weights associated with $Z_{X_{n}}^{RW}[M, {\mathcal O}(k)]$ are
$b_{\Gamma_{(n-k)}}$. It is easiest to start with a generating
function for the insertion of the operators,
\be
Z_{X_{n}}^{RW}[M, \alpha ] = \int D\Phi \, \ex{ i S_{0}} \, \ex{ i
\int_{M} \Sigma_{I} \eta^{I}_{0} - \alpha \frac{d(n)}{2} \eta^{I}_{0} \eps_{IJ}
\eta^{J}_{0} } , \label{gen}
\ee
where $d(n)=(d_{n})^{1/n}$ and 
\be
\Sigma_{I} = \frac{1}{6} \, \Omega_{IJKL}\, \chi^{J}\chi^{K}\chi^{L} -
\frac{1}{2}  \gamma^{AK}_{i}\gamma^{BL}_{j} \eps_{AB} \Omega_{IJKL} \, 
\chi^{J} \f^{j} * d \f^{i} .
\ee
Thus,
\be
Z_{X_{n}}^{RW}[M, {\mathcal O}(k)] = \left. \frac{\partial^{k}}{\partial
\alpha^{k}} Z_{X_{n}}^{RW}[M, \alpha ] \right|_{\alpha =0} .
\ee
One first performs the integration over the modes $\eta^{I}_{0}$ in
(\ref{gen}) to obtain
\be
Z_{X_{n}}^{RW}[M, \alpha ] = (\alpha d(n) )^{n} \int D\Phi \, \ex{i S_{0}} \,
\ex{ \frac{1}{2 d(n) \alpha} \int_{M} \Sigma_{I} \, \eps^{IJ} \, \int_{M}
\Sigma_{J} }
\ee
Consequently,
\bea
Z_{X_{n}}^{RW}[M, {\mathcal O}(k)] & = & \left( | {\mathrm H}_{1}(M,
{\Bbb Z } |  \right)^{n}\frac{k!}{(n-k)!}
\frac{d(n)^{k}}{(2\pi)^{n}} \int_{X_{n}} \sqrt{g} d^{4n} \f_{0} \,
\nonumber \\
& & \;\;\; . \,  \left\langle
\left( \frac{1}{2} \int_{M} \Sigma_{I} \, \eps^{IJ} \, \int_{M}
\Sigma_{J} \right)^{n-k} \right\rangle \label{Xk}
\eea
The integration over the $\chi^{I}$ and $\f^{i}$ fields generates all
the graphs with the appropriate weights. We write,
\be
\left\langle
\left( \frac{1}{2} \int_{M} \Sigma_{I} \, \eps^{IJ} \, \int_{M}
\Sigma_{J} \right)^{k} \right\rangle (X_{n}) = \sum_{\Gamma_{k}}
W_{\Gamma_{k}}(X_{n},\f_{0}) \, I_{\Gamma_{k}}(M) , 
\ee
and clearly
\be
b_{\Gamma_{n}}(X_{n}) = \frac{1}{(2\pi)^{2n}} \int_{X} \sqrt{g} d^{4n}
\f_{0}  \,  W_{\Gamma_{n}}(X_{n},\f_{0})
\,  ,
\ee
as it should. Substituting these expressions back into (\ref{Xk}) one gets
\bea
Z_{X}^{RW}[M, {\mathcal O}(k)]  &= &  \left( | {\mathrm H}_{1}(M, {\Bbb Z}
|  \right)^{n} \frac{k!}{(n-k)!}
\frac{d(n)^{k}}{(2\pi)^{n}} \int_{X_{n}} \sqrt{g} d^{4n} \f_{0} \,
 \nonumber \\
& & \;\;\; . \, 
\sum_{\Gamma_{(n-k)}} 
W_{\Gamma_{(n-k)}}(X_{n},\f_{0}) \, I_{\Gamma_{(n-k)}}(M) . \label{Xk2}
\ee
Now compare this with the
partition function for a hyper-K\"{a}hler manifold $X'_{(n-k)}$,
\bea
Z_{X'_{(n-k)}}^{RW}[M] &=&  \frac{\left( | {\mathrm H}_{1}(M, {\Bbb Z}
| \right)^{(n-k) }}{(2\pi)^{2(n-k)}} \int_{X'_{(n-k)}} \sqrt{g'} d^{4(n-k)}
\f_{0} \, \left\langle
\left( \frac{1}{2} \int_{M} \Sigma_{I} \, \eps^{IJ} \, \int_{M}
\Sigma_{J} \right)^{n-k} \right\rangle \nonumber \\
& = & \frac{\left( | {\mathrm H}_{1}(M, {\Bbb Z} | \right)^{(n-k)
}}{(2\pi)^{2 (n-k)}} \int_{X'_{(n-k)}} \sqrt{g'} d^{4(n-k)}
\f_{0} \nonumber \\
& & \,\, . \, \sum_{\Gamma_{(n-k)}} 
W_{\Gamma_{(n-k)}}(X'_{(n-k)},\f_{0}) \, I_{\Gamma_{(n-k)}}(M). \label{Xn-k}
\eea
We see, therefore, that the same graphs will appear in (\ref{Xk2}) as
in (\ref{Xn-k}). The weights are the same if we define them by
\be
b_{\Gamma_{k}}(X_{n}) = \frac{k!}{(n-k)!}
\frac{d(n)^{k}}{(2\pi)^{n}} \int_{X} \sqrt{g} d^{4n} \f_{0} \,
 W_{\Gamma_{k}}(X_{n},\f_{0}) 
\,  .
\ee
With these observations in hand we have
\be
Z^{RW}_{n}[M, {\mathcal O}(n-k)] = \left( | {\mathrm H}_{1}(M, {\Bbb Z})
| \right )^{n-k} Z_{k}^{RW}[M] .
\ee

\subsection{Equivalence of the RW and LMO Invariants}

One may consider the lambda invariants in the same way as one does the
$Z_{X}^{RW}[M, {\mathcal O}(k)]$, that is one introduces a
$\lambda_{n}^{k}(M)$ which acts on $X$'s. Their definition may be read
off (\ref{ldef}),
\be
\lambda_{n}^{k}(M) = Z_{n}^{RW}[M, {\mathcal O}(n-k)] -
\sum_{m=0}^{k-1} \lambda_{n}^{m}(M) \, Z^{RW}_{n}[S^{3} , {\mathcal
O}(n+m-k)] .
\ee

The aim of this section is to convince the reader that the LMO
invariants $\Omega_{n}\left(M\right)^{(k)}$ (\ref{om}) and the lambda
invariants, $\lambda_{n}^{k}(M)$, are equal
\be
\Omega_{n}\left(M\right)^{(k)} = \lambda_{n}^{k}(M) . \label{lmorw}
\ee
For $b_{1}(M) \geq 1$ we have already seen that the identification is
correct. For $b_{1}(M)=0$, we cannot, at present, completely prove the
equivalence, even at the physical level of rigour. In order to do that
we would have to exhibit surgery formulae for both the Rozansky-Witten
and LMO invariants, formulae which, unfortunately, we do not
know. Instead we offer four good reasons for believing (\ref{lmorw}):

\begin{description}

\item{{\bf Normalization}}

The LMO invariant is designed to be unity for
the 3-sphere. This also motivates the identification of
$\lambda_{n}^{k}(M)$ with $\Omega_{n}(M)^{(k)}$ since, by (\ref{CONSUM}), 
\be
\lambda_{n}^{k}(S^{3})=0. 
\ee
That it is $\lambda_{X}^{k}$ that has this property rather than
$Z_{X}^{RW}$ can be seen in perturbation theory. Since in perturbative
expansions there is a close correspondence between the Rozansky-Witten
invariant and the Chern-Simons invariants there is no reason why the
Rozansky-Witten invariants should vanish for the 3-sphere (see
below for more consequences of this).

\item{{\bf Orientation}}

The $\Omega_{n}(M)^{(k)}$ satisfy\footnote{Beware this is stated
incorrectly in the eprint archive version of Proposition 5.2 in
\cite{LMO}.} \cite{LMO}
\be
\Omega_{n}(^{*}M)^{(k)} = (-1)^{k(b_{1}(M)+1)}\,  \Omega_{n}(M)^{(k)} .
\ee
One can show, inductively, that the $\lambda^{k}_{X}(M)$ behave under
orientation reversal in the same way as $Z_{X}^{RW}[M, {\mathcal
O}(n-k)]$, that is
\be
\lambda^{k}_{n}(M^{*}) = (-1)^{k(b_{1}(M)+1)} \, \lambda^{k}_{n}(M) .
\ee

\item{{\bf Weight Systems}}

By Lemma (4.6) of \cite{LMO} one has
\be
\Omega_{n}\left(M\right)^{(k)} = m^{(n-k)} \, 
\Omega_{k}\left(M\right)^{(k)} , \label{lemma}
\ee
where,
\be
m= \Omega_{1}\left(M\right)^{(0)} = |{\mathrm H}_{1}(M, {\Bbb Z})| .
\ee
For (\ref{lmorw}) to hold one then requires a similar relationship
amongst the $\lambda_{n}^{k}(M)$. For $M$ of rank one or greater
(\ref{lemma}) implies that $\Omega_{n}\left(M\right)^{(k)} =0$ for $k
\neq n$. A similar story holds for $\lambda_{n}^{k}(M)$, namely if $M$
has rank greater than 0, then for $ k \neq n$, $\lambda_{n}^{k}(M)=
0$. We now need to address the rank 0 case. 

Let $M$ be a $\QHS$. Suppose that for all $j \leq n$ and for all $i$ up
to some fixed value $k$ (less than $j$) that
\be
\lambda_{j}\left(M\right)^{(i)} = m^{(j-i)} \, 
\lambda_{i}\left(M\right)^{(i)} . \label{lemma2} 
\ee
We show that this implies that (\ref{lemma2}) holds for $i=k+1$, 

\bea
\lambda_{j}\left(M\right)^{(k+1)} & = & Z_{j}^{RW}[M, {\mathcal O}(j-k-1)] -
\sum_{i=0}^{k} \lambda_{j}^{i}(M) \, Z^{RW}_{j}[S^{3} , {\mathcal
O}(j+i-k-1)] \nonumber \\
& =& m^{j-k-1} \left(Z_{k+1}^{RW}[M] -
\sum_{i=0}^{k} m^{k+1 - i}\lambda_{i}^{i}(M) \, Z^{RW}_{k+1-i}[S^{3}]
\right) \nonumber \\
& = & m^{j-k-1} \left(Z_{k+1}^{RW}[M] -
\sum_{i=0}^{k} \lambda_{i}^{k+1}(M) \, Z^{RW}_{k+1-i}[S^{3}]
\right) \nonumber \\
& = & m^{j-k-1} \lambda_{k+1}^{k+1} .
\eea
To complete the induction we need only notice that $\lambda_{l}^{0}(M) =
Z_{l}^{RW}[M, {\mathcal O}(l)]= m^{l} = m^{l}\lambda_{0}^{0}(M)$ for
all $l$. Consequently, starting with $i=0$, $j=1$ we can prove
inductively that (\ref{lemma2}) holds for all $i$ and for all $j$.

\item{{\bf Connected Sum Properties}}

The LMO invariant \cite{LMO} $\Omega(M)$ satisfies the following
connected sum formula
\be
\Omega_{n}\left(M_{1}\# M_{2}\right)^{(n)} =
\sum_{d_{1}+d_{2}=n}^{\infty} \Omega_{n}\left(M_{1}\right)^{(d_{1})}
\Omega_{n}\left(M_{2}\right)^{(d_{2})} 
\ee
A glance at (\ref{CONSUM}) shows us that the $\lambda_{n}^{n}(M)$
satisfy the same rule under connected sum as the
$\Omega_{n}\left(M\right)^{(n)}$.

\end{description}

\noindent\underline{{\bf The Conjecture}}

The conjectured relationship between the Rozansky-Witten invariants
and LMO can now be stated as
\be
\left| {\mathrm H}_{1}(M, {\Bbb Z})\right|^{n-k} W^{RW}_{X}\left
( Z_{k}^{LMO} (M)  \right)
= \lambda_{X}^{k}(M) \label{conjecture}
\ee

\noindent For $n=k=1$, the left hand side has been shown in \cite{LMO}
to be proportional to the Casson-Lescop-Walker invariant, while
\cite{RW} established that the right hand side is proportional to the
same invariant. The proportionality constants are fixed by the weights
and so the conjectured equality (\ref{conjecture}) has been
established in this case.

\noindent\underline{{\bf Remark (1)}}

\noi The analogy with Chern-Simons theory helps to
understand the connected sum formula and the normalization, since there
we know that \cite{W2}
\be
\frac{Z_{CS}(M_{1}\# M_{2})}{Z_{CS}(S^{3})} =
\frac{Z_{CS}(M_{1})}{Z_{CS}( S^{3})} \, . \, \frac{Z_{CS} \label{csfact}
( M_{2})}{Z_{CS}(S^{3})} .
\ee
With this behaviour we see that the correct ``normalized'' invariants
are $\hat{Z}(M)=Z_{CS}(M)/Z_{CS}( S^{3})$ in Chern-Simons
theory, whence $\hat{Z}(S^{3})=1$. Recall that the $Z^{RW}[M]$ involve
diagrams that arise in Chern-Simons theory, however the $\hat{Z}(M)$
invariants, which satisfy (\ref{csfact}) involve differences of
diagrams for the given manifold with diagrams for $S^{3}$. This is the
nature of the $\lambda_{n}^{k}(M)$ invariants. Denoting the terms
proportional to $1/\kappa^{n}$ in $\hat{Z}$ by $\hat{Z}_{n}$ and in
$Z_{CS} $ by $Z_{CS, \, n}$ we have
an expansion $\hat{Z}_{n}(M) = Z_{CS,\, n}(M) - Z_{CS,\, n-1}(M)Z_{CS,
\, 1}( S^{3}) + \dots $, which is the analogue of (\ref{ldef}).

\noindent\underline{{\bf Remark (2)}}

\noi It is quite impressive that, by this identification, the precursors of
the LMO invariant, the
$\lambda^{k}_{X}$, are coefficients of a vector, in a particular
basis, in the Hilbert space on $S^{2}$.

\section{The RW and Generalized Casson Invariants}\label{gcasson}

Our suggestion for the correct generalization of the $SU(2)$ Casson invariant,
$\lambda_{M}$, to gauge groups $G$ of rank $n$ is
\be
\lambda_{G}(M) = \lambda_{X_{G}}^{n}(M) .
\ee
Part of the motivation for this is of course the relationship between
the $\lambda_{X}^{k}$ and the LMO invariants. Perhaps more importantly
this definition goes some way towards making contact with the work in
\cite{BH} where an $SU(3)$ Casson invariant is defined
rigorously. There the invariant vanishes for $S^{3}$, whereas
$Z_{X_{SU(3)}}^{RW}[S^{3}] \neq 0$. 

The analogy with Chern-Simons
theory allows us to show that $Z_{X_{SU(3)}}^{RW}[S^{3}]$ does not
vanish. Recall that the partition function, $Z^{CS}_{SU(2)}[S^{3},\kappa]$, of
$SU(2)$ Chern-Simons theory for $S^{3}$ is known in closed
form\footnote{The variable $\kappa = k+2$ where $k$ is the level.},
\be
Z^{CS}_{SU(2)}[S^{3},\kappa] = \sqrt{\frac{2}{\kappa}} \, \sin{\left
( \frac{\pi}{\kappa} \right) } . \label{cspart}
\ee
Since the only flat connection on $S^{3}$ is the trivial connection,
the perturbative expansion in $1/\sqrt{\kappa}$ about it should agree with
the large $\kappa$ expansion of (\ref{cspart}). The one-loop, or
Ray-Singer Torsion contribution goes like $\kappa^{-3/2}$, so that the loop
expansion is
\be
Z^{CS}_{SU(2)}[S^{3},\kappa] = \frac{\sqrt{2} \pi}{\kappa^{3/2}} \left( 1 +
\sum_{n=1}^{\infty} \frac{a_{n}}{\kappa^{n}} \right) 
\ee
and a comparison with (\ref{cspart}) shows that $a_{1}=0$ while $a_{2}
\neq 0$. The vanishing of $a_{1}$ tells us that there is no
contribution from the $\Theta$ diagram, while the non-vanishing of
$a_{2}$ tells us that the sum of the double theta plus the Mercedes
Benz diagram does not vanish. We know then that the $I_{\Gamma}$
associated with these diagrams do not vanish. As far as the group
theory factor of these
diagrams is concerned, it is proportional to the square of the
quadratic casimir 
of the group, which also does not vanish for $SU(n)$. Consequently,
these diagrams will contribute in the case of $SU(3)$.

Why does the path integral that is `designed' to yield the generalized
Casson invariant not do so? The Casson invariant, from the
gauge theory point of view of \cite{Ta} or \cite{BH}, is such that
the trivial connection is always `excised' when it comes to
performing a count of (perturbed) flat connections on
$\ZHS$'s. This is the reason that the Casson invariant and its
generalization vanish for $S^{3}$. There is, however, no such
directive in the path integral for the supersymmetric gauge theory
that was analyzed in \cite{SW} and consequently, no such directive in
the path integral formulation of \cite{RW}. It is
clear, then, that it is necessary to subtract off the contribution, if
any, of the trivial connection in order to arrive at the generalized
Casson invariant. 

How are we to perform the excision? The perturbations, both in
\cite{Ta} and \cite{BH} for $\ZHS$'s, are designed so that the product
connection is isolated from the other (perturbed) flat connections that
contribute to the invariant. This is also the case for the trivial
representation in Walker's definition \cite{Wa} of the Casson invariant for
$\QHS$'s. Since the trivial connection contribution in these cases
can be isolated its contribution can be subtracted if it is
known. Unfortunately, the problem is that we do not know it. If the
formula,
\be
\lambda_{X}^{n}(M) = Z_{X}^{RW}[M] - \sum_{k=0}^{n-1}
\lambda_{X}^{k}(M)Z_{X}^{RW}[ S^{3} , {\mathcal O}(k)] 
\ee
is the correct definition for the Casson invariant, this suggests that the
contribution of the trivial connection is `universal' in that,
regardless of $M$ one is subtracting out the contribution of the
trivial connection in $S^{3}$ (up to factors that depend on the
cohomology of the $\QHS$).

If $M$ is not a $\QHS$ then the product connection is not isolated
from the rest of the moduli space of flat connections and one cannot
``cleanly'' subtract off its contribution to the invariant. However,
for such an $M$, we have 
\be
\lambda_{X}^{n}(M) = Z_{X}^{RW}[M] ,
\ee
since $Z_{X}^{RW}[M, {\mathcal O}(n-k)] = 0$, $\forall k \neq n$ (since
for $b_{1}(M) >0$ the insertion of these observables gives zero).

It is difficult to completely fix the relationship between the invariant of
\cite{BH} and the one that we are proposing for the $SU(3)$ Casson
invariant. However, we suspect they are closely related to each
other. The reason we believe
this is that both are generalizations of the same object in the
$SU(2)$ case. \cite{BH} generalize the gauge theoretic construction of
Taubes \cite{Ta}. However, Taubes's approach is equivalent to the
physics approach in the $SU(2)$ case as presented in \cite{W1} and
\cite{BT1}. Here one can clearly excise the trivial connection and give
a treatment more in line with the mathematical one. On the other
hand, Rozansky and Witten are generalizing an alternative physics approach
to the $SU(2)$ invariant. With the small caveat made above, they are
evaluating an equivalent path integral to the one proposed in \cite{BT1}. 

As a small check we note that both of the invariants
vanish on $S^{3}$ and they are both insensitive to the orientation of $M$.

\noindent{\bf Note Added (August):} Boden and Herald have shown that their
invariant, which we denote by $\lambda_{{\mathrm BH}}$,
satisfies\footnote{We would like to
thank Chris Herald for informing us of this prior to publication.},
\cite{BH2},
\be
\lambda_{{\mathrm BH}}(M_{1} \# M_{2}) = \lambda_{{\mathrm BH}}(M_{1})+
\lambda_{{\mathrm BH}}(M_{2}) + 4
\lambda_{SU(2)}(M_{1})\lambda_{SU(2)}(M_{2}) ,
\ee
where $\lambda_{SU(2)}$ is normalized as in \cite{Wa}. While this is
consistent with the invariants $\lambda_{X_{SU(3)}}$ and
$\lambda_{{\mathrm BH}} $ being proportional
\be
\lambda_{X_{SU(3)}}(M) = \a\, \lambda_{{\mathrm BH}}(M), \label{rel}
\ee
and
\be
\lambda_{X_{SU(3)}}^{1}(M) = 2\sqrt{\a} \; \lambda_{SU(2)}(M) ,
\ee
for some $\a$, we do not believe that they can be related in such a
simple way. Recently, it has been
shown that the Boden-Herald invariant is not of finite type of degree
$\leq 6$, \cite{BHKK} (their theorem 6.16) and so if (\ref{rel}) were
true then our conjecture would be false. 

\noindent{\bf Note Added (August):} By making use of supersymmetry
and other physics inspired arguments Paban, Sethi and Stern \cite{PSS}
have determined the integral of the Euler density over the $SU(2)$
n-monopole moduli space and have found that it is equal to
$n$. However, integrals of other densities are still not known.

\section{The Appendices}
\appendix
\section{Some Properties of Hyper-K\"{a}hler Manifolds}\label{hkprops}

Generically the holonomy group of a real $m$ dimensional Riemannian manifold is
${\mathrm SO}(m)$. If the manifold is complex of complex dimension $p$
and the metric is hermitian then the holonomy lies in ${\mathrm U}(p)
\subset  {\mathrm SO}(2p)$.
If in addition $X$ is a hyper-K\"{a}hler manifold ($\dim_{{\Bbb R}}X=4n$) then
there is a hermitian metric such that the Levi-Civita connection lies
in the ${\mathrm Sp}(n)$ subgroup of ${\mathrm U}(2n) \subset {\mathrm
SO}(4n)$. The
complexified tangent bundle decomposes as
\be
TX_{{\Bbb C}} = TX \otimes_{{\Bbb R}} {\Bbb C} = V \otimes_{{\Bbb C}} S ,
\label{decomp}
\ee
where $V$ is a rank $2n$ complex vector bundle with structure group
${\mathrm Sp}(n)$ and $S$ is a trivial rank 2 complex vector bundle
with structure group ${\mathrm Sp}(1)$. The Levi-Civita connection is
a connection on $V$ and is the trivial connection on $S$. ${\mathrm
Sp}(1)$ labels are A, B, $\dots$, $= 1,2$, and there is an invariant
antisymmetric tensor $\eps_{AB}$ with inverse $\eps^{AB}$,
\be
\eps^{AC}\eps_{CB} = \d^{A}_{B} .
\ee
${\mathrm Sp}(n)$ labels are I, J, $\dots$ $= 1, \dots, 2n$, and there
is also an
invariant antisymmetric tensor $\eps_{IJ}$ with inverse $\eps^{IJ}$,
\be
\eps^{IK}\eps_{KJ}=\d^{I}_{J} .
\ee

Local coordinates on $X$ will be denoted $\f^{i}$ and the Riemannian
metric is $g_{ij}$. The fact that the tangent bundle decomposes as in
(\ref{decomp}) means that there exist covariantly constant tensors
$\gg^{AI}_{i}$ and $\gg_{AI}^{i}$ that describe the maps from $V
\otimes S$ to $TX_{{\Bbb C}}$ and vice versa,
\bea
& & \gg_{i}^{AI} : V \otimes S \rightarrow TX_{{\Bbb C}} \nonumber \\
& & \gg_{AI}^{i} : TX_{{\Bbb C}} \rightarrow V \otimes S .
\eea
These maps are inverses in the sense that
\be
\gg_{i}^{AI}\gg^{i}_{BJ} = \d^{A}_{B} \d^{I}_{J} .
\ee

Using these tensors one
may express the Riemann curvature tensor as
\be
R_{ijkl} = - \gg^{AI}_{i} \gg^{BJ}_{j} \gg^{CK}_{k} \gg^{DL}_{l}
\eps_{AB} \eps_{CD} \Omega_{IJKL} , \label{iso}
\ee
where $\Omega_{IJKL}$ is completely symmetric in the indices. A useful
relationship is
\be
g_{ij}\gg^{j}_{AI} = \eps_{AB}\eps_{IJ}\gg_{i}^{BJ} .
\ee

Fix on a complex structure so that $\f^{I}$ are holomorphic
coordinates on $X$ with respect to this complex structure. Then we may
take
\be
\gg_{A J}^{I} = \d_{A1} \d^{I}_{J}, \;\;\;\; \gg_{A I}^{\overline{I}}
= \d_{A 2} g^{\overline{I}J}\eps_{JI}.
\ee
In such a preferred complex structure, the tensor
\be
T^{\overline{J}}_{\; J} = g^{\overline{J}K}\eps_{KJ}
\ee
maps $T^{(1,0)}X$ to $T^{(0,1)}X$ while
\be
T_{\; \overline{J}}^{ J} = \eps^{JK}g_{K\overline{J}},
\ee
maps $T^{(0,1)}X$ to $T^{(1,0)}X$. Since the $T$ tensors are inverses
of each other,
\be
T^{\overline{J}}_{\; J} \, T^{J}_{\; \overline{K}} =
\d^{\overline{J}}_{\overline{K}} , \;\;\; {\mathrm and}, \;\;\;
T^{J}_{\; \overline{J}}T^{\overline{J}}_{\; K} = \d^{J}_{K}, 
\ee
they provide an isomorphism between $T^{(1,0)}X$ and $T^{(0,1)}X$.

The holomorphic symplectic two form $\eps$ is covariantly constant
with respect to the Levi-Civita connection on $V$,
\be
\partial_{K}\eps_{IJ} - \Gamma^{L}_{KI}\eps_{LJ}-
\Gamma^{L}_{KJ}\eps_{IL}= 0. \label{covcon}
\ee
In the preferred complex structure one finds that
\be
\Omega_{IJKL} = - R_{I\overline{J}K\overline{L}}\, T^{\overline{J}}_{J}
\, T^{\overline{L}}_{L} .\label{or}
\ee
At this point it is not completely transparent that $\Omega_{IJKL}$ is
totally symmetric in the labels, however, it is indeed so. We pause to
prove this. Since $\eps_{IJ}$ is holomorphic we have
\be
\partial_{\overline{K}}\partial_{K}\eps_{IJ}=
\partial_{K}\partial_{\overline{K}}\eps_{IJ}= 0.
\ee
~From (\ref{covcon}), this means that
\be
\partial_{\overline{K}}\Gamma^{L}_{KI}\eps_{LJ} =
\partial_{\overline{K}} \Gamma^{L}_{KJ}\eps_{LI},
\ee
however for a K\"{a}hler manifold one has
\be
R^{I}_{\; JK \overline{L}} = \partial_{\overline{L}} \Gamma^{I}_{JK} ,
\ee
so that we have shown
\be
\eps_{IL}R^{L}_{\; JK \overline{L}} = \eps_{KL}R^{L}_{\; JI
\overline{L}} .
\ee
Hence $\eps_{IL}R^{L}_{\; JK \overline{L}}$ is symmetric in $I$, $J$
and $K$. On the other hand, from (\ref{iso}) and (\ref{or}), we have
\be
\Omega_{IJKL}= \eps_{NJ}R^{N}_{IK \overline{L}} \,
T^{\overline{L}}_{L},
\ee
which shows that $\Omega_{IJKL}$ is totally symmetric.

\section{Berezian Integration}\label{sec.manipulations}

Let $V$ be a vector space. By a polynomial (bosonic) function on $V$,
we mean an element of $S(V^*)$, the symmetric tensor algebra of  $V^*$, the
dual of $V$. By a Grassmann (fermionic) function on $V$, we mean an
element of $\Lambda(V^*)$, the exterior algebra.

The Berezian integral of an element of $\Lambda(V^*)$, is its
projection to the top dimensional piece (provided $V$ is finite
dimensional).  It is a number, provided we have a metric and
orientation of $V$ (which yields a metric and
orientation, and hence a trivialization, of the top exterior power).

The rules for Berezin integration mean that the non-zero
linear map $T:\Lambda V^{*} \rightarrow {\Bbb R}$ is indeed an `integral'
and will be denoted as $\int d^{m}\theta$. Here, the $\theta^{\mu}$
form a basis of $\Lambda^{1}(V^*)$. We normalize this integral
in the following way, we let
\be
\int d^{m}\theta \, \theta^{\mu_{1}} \dots \theta^{\mu_{m}} =
\varepsilon^{\mu_{1} \dots \mu_{m}} \label{brez}
\ee
where
\be
\varepsilon = \left\{ \begin{array}{cl}
1 & {\mathrm even \,\, permutation}\\
-1 & {\mathrm odd \,\, permutation} \\
0 & {\mathrm otherwise}
\end{array} \right.
\ee
Given $A \in \Lambda^{2}V$, $A= -1/2 \, A_{ij}e^{i}\wedge e^{j}$, the
Pfaffian of $A$ is defined to be the number,
\be
{\mathrm Pfaff}(A) = \frac{(-1)^{m}}{2^{m}m!} \varepsilon^{i_{1} \dots
i_{m}} A_{i_{1}i_{2}} \dots A_{i_{2m-1}i_{2m}}.
\ee
By making use of the Berezin integral, one may also write this as,
\bea
{\mathrm Pfaff}(A) &=& T \left( \ex{A}\right) \nonumber \\
& = & \int d^{m}\theta \, \ex{-\frac{1}{2} \theta^{i} A_{ij}
\theta^{j}} .
\eea

The exterior algebra satisfies, for the top wedge product,
\bea
dx^{\mu_{1}} \wedge \dots \wedge dx^{\mu_{m}}
&= & \varepsilon^{\mu_{1}
\dots \mu_{m}}dx^{1}\wedge \dots \wedge dx^{m} ,\nonumber \\
&= & \int d^{m} \theta \, \theta^{\mu_{1}} \dots \theta^{\mu_{m}} \;
d^{m}x \label{vol}
\eea
the second equality follows from (\ref{brez}) and we have chosen the
orientation
\be
d^{m}x = dx^{1} \dots dx^{m}.
\ee
Consequently, for any
top form
\be
f= f_{\mu_{1} \dots \mu_{m}}\, dx^{\mu_{1}}
\dots dx^{\mu_{m}},
\ee
we have
\be
f = \left( \int d^{m}\theta \, f_{\mu_{1} \dots \mu_{m}}\,  \theta^{\mu_{1}}
\dots \theta^{\mu_{m}} \right) \, d^{m} x . \label{topf}
\ee

So far we have not included the notion of a metric. As it stands in
Riemannian geometry $\varepsilon$ is a density and not a tensor. Fix
on a metric $g_{\mu \nu}$. Now one sees that $\varepsilon/\sqrt{g}$,
where $ g = \det{ g_{\mu \nu}}$, is a tensor. We introduce a new
measure for the Berezin integration
\be
f = \int d\mu_{g}(\theta) \; f_{\mu_{1} \dots \mu_{m}} \, \theta^{\mu_{1}}
\dots \theta^{\mu_{m}} \, \sqrt{g} \,d^{m} x , \label{topfmet}
\ee
so that $d\mu_{g}(\theta) = d^{m}\theta/\sqrt{g}$, or put another way
\be
\int d\mu_{g}(\theta) \, \theta^{\mu_{1}} \dots \theta^{\mu_{m}} =
\frac{\varepsilon^{\mu_{1} \dots \mu_{m}}}{\sqrt{g}} . \label{brezmet}
\ee

The more structures that are introduced on the manifold the more
variations that are available on this theme. The first refinement is
to consider manifolds of dimension $2m$ that come equipped with a
complex structure. In
such a situation we can refine the formula (\ref{topf}) for forms of
degree $(m,k)$
\be
g = g_{I_{1} \dots I_{m}\, \bar{J}_{1} \dots \bar{J}_{k}}\,  
dz^{I_{1}}\dots dz^{I_{m}} \, d\bar{z}^{\bar{J}_{1}} \dots
d\bar{z}^{\bar{J}_{k}} , 
\ee
to
\be
g = \left( \int d^{m}\theta \, g_{I_{1} \dots I_{m}\, \bar{J}_{1}
\dots \bar{J}_{k}}\,  \theta^{I_{1}}
\dots \theta^{I_{m}} \right) \, d^{m}z \, d\bar{z}^{\bar{J}_{1}} \dots
d\bar{z}^{\bar{J}_{k}} . \label{brezc}
\ee

Our interest is in hyper-K\"{a}hler manifolds of real dimension $4n$,
where we are assured of
the existence of a holomorphic symplectic 2-form $\eps$
\be
\eps = - \frac{1}{2}\eps_{IJ}dz^{I} dz^{J}
\ee
which is non-degenerate, $\eps^{n} \neq 0$. The inverse matrix
$\eps^{IJ}$ is defined by $\eps^{IK}\eps_{KJ} = \d^{I}_{J}$. The
analogue of (\ref{vol}) is
\bea
dz^{I_{1}} \dots dz^{I_{2n}} &=& \varepsilon^{I_{1} \dots I_{2n}} \,
dz^{1} \dots dz^{2n} \nonumber \\ 
& =& \left( \int d^{2n}\eta \, \eta^{I_{1}} \dots \eta^{I_{2n}}
\right) \, d^{2n}z .
\eea
Furthermore we have the definition
\be
\eps^{I_{1} \dots I_{2n}} = \varepsilon^{I_{1} \dots I_{2n}}\,
\frac{{\mathrm Pfaff}(\eps)}{\det{(\eps)}}
\ee
so that
\bea
dz^{I_{1}} \dots dz^{I_{2n}} &=& \eps^{I_{1} \dots I_{2n}}\, {\mathrm
Pfaff}(\eps)\, dz^{1} \dots dz^{2n} \nonumber \\
& =& \eps^{I_{1} \dots I_{2n}} \, \eps^{n} .
\eea
Thus the analogue of (\ref{brezc}) for $(2n,0)$ forms is
\be
f = \left( \int d^{2n}\eta \, f_{I_{1} \dots I_{2n}} \, \eta^{I_{1}}
\dots \eta^{I_{2n}} \right) \, d^{2n}z ,
\ee
but, since we have the holomorphic symplectic 2-form at our
disposal, we have the analogue of the metric dependent measure
(\ref{topfmet})
\be
f=  \left( \int d\mu(\eta) \, f_{I_{1} \dots I_{2n}} \, \eta^{I_{1}}
\dots \eta^{I_{2n}} \right) \, \eps^{n}.
\ee
The new measure (we do not exhibit the $\eps$ dependence in the
measure $d\mu$ which therefore should be, more correctly, denoted by
$d\mu_{\eps}$) is
\be
d\mu(\eta) = d^{2n}\eta \, \frac{{\mathrm Pfaff}(\eps)}{\det{(\eps)}} .
\ee

The isomorphism between $TX^{(1,0)}$ and $TX^{(0,1)}$ means that we
also have  the following, anti-holomorphic, 2-form
\be
\bar{\eps} = - \frac{1}{2} \, \bar{\eps}_{\bar{I}\bar{J}} \,
d\bar{z}^{\bar{I}}\, d\bar{z}^{\bar{J}} = -\frac{1}{2} \eps_{IJ} \,
T^{I}_{\bar{I}} T^{J}_{\bar{J}}\, d\bar{z}^{\bar{I}}\, d\bar{z}^{\bar{J}} .
\ee 
Hence,
\bea
d\bar{z}^{\bar{I}_{1}} \dots d\bar{z}^{\bar{I}_{2n}} &=&
\bar{\eps}^{\bar{I}_{1} \dots \bar{I}_{2n}}\, {\mathrm
Pfaff}(\bar{\eps})\, d\bar{z}^{1} \dots d\bar{z}^{2n} \nonumber \\
& =& \bar{\eps}^{\bar{I}_{1} \dots \bar{I}_{2n}} \, \bar{\eps}^{n} .
\eea

\subsection{Normalization of Zero Modes}

The normalization of the path integral measure for zero modes that we
adopt is best stated in the following manner. For each Grassmann
valued section of $\f^{*}_{0}V$, denoted by $\eta^{I}$ we demand that
\be
\int d\mu(\eta )\, \ex{-\frac{1}{2} \eta^{I}\eps_{IJ}\eta^{J}} = 1 .
\ee
This is in contrast to the more standard measure, that is used above,
for which the integral over the zero modes is normalized as
\be
\int d^{2n} \eta \, \ex{-\frac{1}{2} \eta^{I}\eps_{IJ}\eta^{J}} =
{\mathrm Pfaff}(\eps ) ,
\ee
and the relationship between the two is clearly
\be
d^{2n}\eta = {\mathrm Pfaff}(\eps) d\mu(\eta).
\ee
One important property that we will make use of is a change of
variables formula
\be
\int d\mu(\eta) \, f(\eps_{IJ}\eta^{J}) = {\mathrm Pfaff}(\eps) \int
d^{2n} \eta \, f(\eta_{I}) .
\ee

\subsection{The Euler Class and Grassmann Integration}
Here we briefly review the construction for
expressing the Euler characteristic of a compact closed manifold in
a form which involves Grassmann variables and Grassmann integration
and which is suitable to our needs. We recall that the Euler class
${\mathbf e}(E)$ of a real vector bundle $E\rightarrow X$ of rank
$2m$, with a given connection $A$ whose curvature is $F_{A} = dA +
A^{2}$, is defined to be the cohomology class
\be
{\mathbf e}(E) = \frac{1}{(2\pi)^{n}}{\mathrm Pfaff}(F_{A})
. \label{eclass}
\ee
The Euler characteristic, $\chi(E)$, of $E$ is
\be
\chi(E) = \int_{X} \, {\mathbf e}(E) .
\ee
When $E$ is the (real) tangent bundle, $TX$, of $X$, we will write
Euler class as
${\mathbf e}(TX)$ and the Euler characteristic as ${\mathbf e}(X)=
\int_{X}{\mathbf e}(TX)$. Using the
rules of Grassmann integration we see that the Euler class may be
represented as
\be
{\mathbf e}(E) = \frac{1}{(2\pi)^{m}}\int  d^{2m}\eta \;
\ex{-\frac{1}{2} \eta_{a}F^{ab}\eta_{b} } .
\ee

\noi {\bf Claim}\label{c1}: The Euler Class for a compact
closed Riemannian manifold $X$ of dimension $2m$ is
\be
{\mathbf e}(TX) = \frac{1}{\sqrt{g}} \left(\int d^{2m}\chi\,
d^{2m}\psi \; \ex{\frac{1}{4} R_{\mu \nu \kappa \lambda} \chi^{\mu}
\chi^{\nu} \psi^{\kappa}\psi^{\lambda}}\right)\, d^{2m}x . \label{ecr}
\ee

\noi {\bf proof:}

The curvature two-form for the tangent bundle is
\bea
F^{ab}& = & R^{ab} \nonumber \\
      &=& \frac{1}{2} R^{ab}_{\;\;\; \mu \nu }dx^{\mu}dx^{\nu}
\nonumber \\
& = & \frac{1}{2} e^{a \sigma}e^{b \rho} \, R_{\sigma \rho \mu \nu}
dx^{\mu} dx^{\nu} \, ,
\eea
where $e^{a} = e^{a}_{\mu}dx^{\mu}$ is a section of the orthonormal
frame bundle. By (\ref{topf}) we have that top form part of the
exponential satisfies
\be
\ex{\frac{1}{2}\eta_{a}R^{ab}\eta_{b}} = \int d^{2m}\psi \;
\ex{\frac{1}{4}R^{ab}_{\;\;\; \mu \nu}\eta_{a}\eta_{b}
\psi^{\mu}\psi^{\nu}} \, d^{2m}x .
\ee
Thus,
\bea
{\mathbf e}(TX)& =& \frac{1}{(2\pi)^{m}} \left(\int d^{2m}\eta \,
d^{2m}\psi\; \ex{\frac{1}{4}R^{ab}_{\;\;\; \mu \nu}\eta_{a}\eta_{b}
\psi^{\mu}\psi^{\nu}} \right)\, d^{2m}x \nonumber \\
&=&  \frac{1}{(2\pi)^{m}} \left(\int d^{2m}\eta \,
d^{2m}\psi \; \ex{\frac{1}{4}R_{\sigma \rho \mu \nu}
e^{a\sigma}\eta_{a} e^{b\rho}\eta_{b}
\psi^{\mu}\psi^{\nu}}\right) \, d^{2m}x . \label{eXep}
\eea
Let
\be
\eta_{a} = e_{a \mu} \chi^{\mu} , \label{etoc}
\ee
and the Jacobian for such a change of variables is
\bea
d^{2m}\eta & = & \det{(e_{a \mu})}^{-1} \, d^{2m} \chi
\nonumber \\
& = & \frac{1}{\sqrt{g}} \, d^{2m} \chi . \label{jac}
\eea
Making the change of variables (\ref{etoc}) in (\ref{eXep}) and
keeping in mind the Jacobian (\ref{jac}) proves the claim.

\noi {\bf Claim:} The Euler Class of a compact closed
hyper-K\"{a}hler manifold $X$ of real dimension $4n$ is
\be
{\mathbf e}(TX) = \sqrt{g} \left( \int d\mu (\chi_{\alpha})\;  \ex
{ \frac{1}{24}\Omega_{IJKL}\chi^{I}_{\alpha}\, \chi^{J}_{\beta}\,
\chi^{K}_{\gamma}\, \chi^{L}_{\d}\, \eps^{\alpha \beta \gamma \d} } \right)
d^{4n}x. \label{c2}
\ee
\noi {\bf Proof:}

The manifold in question carries a hyper-K\"{a}hler
structure and so there is a refinement for the Riemann curvature
tensor that was explained previously. Fix on a prefered
complex structure for $X$. Let $x^{\mu}$ be local coordinates and in
the prefered complex structure let $x^{I}$ be the holomorphic and
$x^{\overline{I}}$ be the anti-holomorphic coordinates. The Riemann
curvature tensor, $R_{\mu \nu \rho \sigma}$, vanishes unless the pairs
of indices $(\mu , \nu)$ and $(\rho , \sigma)$ are of $(I,
\overline{J})$ or $(\overline{I}, J)$ type. Consequently,
\be
\frac{1}{4}R_{\mu \nu \kappa \lambda}\, \chi^{\mu}\chi^{\nu}
\psi^{\kappa} \psi^{\lambda} = R_{I\overline{J}K\overline{L}}\, \chi^{I}
\chi^{\overline{J}} \psi^{K} \psi^{\overline{L}} .
\ee
Hence, the Euler class may now be expressed as
\be
{\mathbf e}(TX) = \frac{1}{\sqrt{g}} \left(\int d^{4n}\chi\,
d^{4n}\psi \; \exp{\left(R_{I \overline{J}K \overline{L}} \chi^{I}
\chi^{\overline{J}} \psi^{K}\psi^{\overline{L}}\right)}\right)\, d^{4n}x
. \label{eci}
\ee
Since $T_{\; I}^{\overline{J}}$ provides an isomorphism between
$T^{(0,1)}X$ and $T^{(1,0)}X$ we may change variables and let
\bea
\chi^{\overline{J}} & =& T^{\overline{J}}_{\; J}\, \chi^{J}_{2}\nonumber \\
\psi^{\overline{L}} & =& T^{\overline{L}}_{\; L}\, \chi^{L}_{4} .
\eea
The measures go to
\bea
d^{2n}\chi^{\overline{J}} &=& \det{\left(T^{\overline{J}}_{\;
J}\right)}^{-1}\, d^{2n}\chi_{2}^{J} \nonumber \\
d^{2n}\psi^{\overline{L}} &=& \det{\left(T^{\overline{L}}_{\;
L}\right)}^{-1}\,  d^{2n}\chi_{4}^{L}.
\eea
However,
\be
\det{\left(T^{\overline{J}}_{\; J}\right)}^{-2} = \det{(g_{\mu \nu})} \;
{\mathrm Pfaff}(\eps_{IJ})^{4} ,
\ee
so that the measure in (\ref{eci}) becomes
\be
\frac{1}{\sqrt{g}}d^{4n}\chi\, d^{4n}\psi = \sqrt{g} \,
\prod_{\alpha =1}^{4}d\mu(\chi_{\alpha}),
\ee
where we have relabeled the fields as $\chi^{I}= \chi^{I}_{1}$ and
$\psi^{I}= \chi^{I}_{3}$. The exponent in (\ref{eci}) is
\be
R_{I\overline{J}K\overline{L}}\, \chi^{I}_{1} \, T^{\overline{J}}_{\;
J}\, \chi^{J}_{2} \, \chi^{K}_{3}\,  T^{\overline{L}}_{\; L}\, \chi^{L}_{4} .
\ee
But,
\be
\Omega_{IJKL} = R_{I\overline{J}K\overline{L}}\, 
T^{\overline{J}}_{\;  J} \, T^{\overline{L}}_{\; L},
\ee
so we are done. 

\noi Note that to prove (\ref{c2}) we do not need to pick a prefered
complex structure. It was expedient to do so here as this is the way
the objects arise in the text.

\noindent \bf{Nathan Habegger}:\\
UMR 6629 du CNRS, Universit\'e de Nantes \\
D\'epartement de Math\'e\-matiques \\ 2 rue de la Houssini\`ere \\ BP
92208\\44322 NANTES Cedex, France  \\
\bf{habegger\char'100math.univ-nantes.fr}

\noindent \bf{George Thompson}:\\
Abdus Salam International Centre for Theoretical Physics\\
P.O.Box 586\\
TRIESTE 34100\\
Italy\\
\bf{thompson\char'100ictp.trieste.it}


\end{document}